\documentclass[12pt]{article}
\usepackage{a4, amsfonts, amssymb, latexsym, amscd, amsmath}
\begin{document}
\newtheorem{theorem}{Theorem}[section] 
\newtheorem{prop}[theorem]{Proposition}
\newtheorem{lemma}[theorem]{Lemma}
\newtheorem{corollary}[theorem]{Corollary}
\newtheorem{conj}[theorem]{Conjecture}
\newcommand{\pf}{{\bf Proof: }}
\newcommand{\finpf}{\hfill$\Box$}
\newcommand{\Jac}{\mathrm{Jac}}
\newcommand{\defn}{{\bf Definition: }}
\newcommand{\tabeq}{\!\!\!\!&=&\!\!\!\!}
\newcommand{\tabneq}{\!\!\!\!&\neq&\!\!\!\!}
\newcommand{\tabplus}{\!\!\!\!&+&\!\!\!\!}
\newcommand{\tabminus}{\!\!\!\!&-&\!\!\!\!}
\newcommand{\tabspace}{\!\!\!\!& &\!\!\!\!}
\newcommand{\tabtimes}{\!\!\!\!&\times&\!\!\!\!}
\newcommand{\tabcolon}{\!\!\!\!&\colon&\!\!\!\!}
\newcommand{\tabequiv}{\!\!\!\!&\equiv&\!\!\!\!}
\newcommand{\tabdivide}{\!\!\!\!&\div&\!\!\!\!}
\newcommand{\tabimplies}{\!\!\!\!&\Rightarrow&\!\!\!\!}
\newcommand{\tabiff}{\!\!\!\!&\Leftrightarrow&\!\!\!\!}
\newcommand{\tabmapsto}{\!\!\!\!&\mapsto&\!\!\!\!}
\newcommand{\tabnotin}{\!\!\!\!&\notin&\!\!\!\!}
\newcommand{\tabcomma}{\!\!\!\!&,&\!\!\!\!}  
\newcommand{\tabin}{\!\!\!\!&\in&\!\!\!\!}
\newcommand{\mcA}{\mathcal A}
\newcommand{\mcC}{\mathcal C}
\newcommand{\mcF}{\mathcal F}
\newcommand{\mcD}{\mathcal D}
\newcommand{\mcG}{\mathcal G}
\newcommand{\mcL}{\mathcal L}
\newcommand{\mcO}{\mathcal O}
\newcommand{\mco}{\mathcal O}
\newcommand{\mcp}{\mathcal P}
\newcommand{\mka}{\mathfrak a}
\newcommand{\mkp}{\mathfrak p}
\newcommand{\mko}{\mathfrak o}
\newcommand{\mkD}{\mathfrak D}
\newcommand{\End}{\mathrm{End}}
\newcommand{\Aut}{\mathrm{Aut}}
\newcommand{\Hom}{\mathrm{Hom}}
\newcommand{\Div}{\mathrm{Div}}
\newcommand{\Pic}{\mathrm{Pic}}
\newcommand{\ddiv}{\mathrm{div}}
\newcommand{\Gal}{\mathrm{Gal}}
\newcommand{\tr}{\mathrm{tr}}
\newcommand{\p}{\mathbb P}
\newcommand{\Q}{\mathbb Q}
\newcommand{\Qbar}{\overline{\mathbb Q}}
\newcommand{\Z}{\mathbb Z}
\newcommand{\F}{\mathbb F}
\newcommand{\h}{\mathbb H}
\newcommand{\C}{\mathbb C}
\newcommand{\R}{\mathbb R}
\newcommand{\T}{\mathbb T}
\newcommand{\hl}{\hline}
\newcommand{\vl}{\vline}
\newcommand{\lb}{\left}
\newcommand{\rb}{\right}
\newcommand{\im}{\mathrm{im}\:}
\newcommand{\spec}{\mathrm{spec}}
\newcommand{\disc}{\mathrm{disc}}
\newcommand{\rem}{{\bf Remark: }}
\newcommand{\GL}{\mathrm{GL}}
\newcommand{\SL}{\mathrm{SL}}
\newcommand{\chr}{\mathrm{char}}
\newcommand{\Supp}{\mathrm{Supp}}     
\renewcommand{\Sp}{\mathrm{Sp}}

\begin{center}
\addtocounter{footnote}{1}
\addtocounter{footnote}{1}
\large{\bf{CURVES OF GENUS 2 WITH $\sqrt2$ MULTIPLICATION}}\\
\vspace{0.25in}
\normalsize{}
PETER R. BENDING\\
\vspace{0.25in}
\end{center}  
\begin{abstract}
We give a universal family of curves of genus 2 whose jacobians have $\sqrt2$ multiplication fixed by the Rosati involution, and several results based on it. 
\end{abstract}            
\section{Introduction}\label{intro}
Let $F$ be a field of degree $g$ over $\Q$, let $N$ be a positive integer, and let $f$ be a newform on $\Gamma_1(N)$ whose field of Fourier coefficients is $F$ ($F$ is a totally real field if $f$ is on $\Gamma_0(N)$, and is a CM field otherwise). Then, the abelian variety associated via the Shimura construction to $f$ is of dimension $g$, defined over $\Q$, and its ring of endomorphisms defined over $\Q$ is isomorphic to an order in $F$ (in particular, it is simple over $\Q$). It is conjectured that an abelian variety which is of dimension $g$, defined over $\Q$, simple over $\Q$, and whose ring of endomorphisms defined over $\Q$ contains a subring which is isomorphic to an order in $F$, is isogenous over $\Q$ to such an abelian variety. Throughout this paper, we will refer to this conjecture as the generalised Shimura-Taniyama-Weil conjecture. 

A. Wiles and R. Taylor have proved the generalised Shimura-Taniyama-Weil conjecture for semistable elliptic curves, implying Fermat's Last Theorem by work due to K. Ribet. F. Diamond has proved it for elliptic curves which are semistable at $3$ and $5$. Very recently, B. Conrad, F. Diamond and R. Taylor have proved it for elliptic curves whose conductor is not divisible by $3^3$. More generally, N. I. Shepherd-Barron and R. Taylor have proved it for abelian varieties such that the ring of integers of $F$ contains a prime whose norm is $3$ or $5$, under some technical conditions (see~\cite{ST}, Theorem 4.2).

The simplest abelian varieties apart from elliptic curves are abelian surfaces, in particular jacobians of curves of genus $2$. If $\alpha$ is an algebraic integer, we say that an abelian variety {\em has $\alpha$ multiplication} if it has an endomorphism which is killed by the minimum polynomial of $\alpha$. This paper is a summary of the author's doctoral thesis, which looks at curves of genus $2$ whose jacobians have $\sqrt2$ multiplication, with the generalised Shimura-Taniyama-Weil conjecture in mind. Unfortunately, N. I. Shepherd-Barron and R. Taylor's theorem cannot be used when $F$ is $\Q(\sqrt2)$, since the ring of integers of $\Q(\sqrt2)$ does not contain a prime whose norm is $3$ or $5$.

We give a universal family of curves of genus $2$ whose jacobians have $\sqrt2$ multiplication fixed by the Rosati involution, and use it to compile a table of curves of genus $2$ whose jacobians satisfy the conditions of the generalised Shimura-Taniyama-Weil conjecture. For each newform $f$ on $\Gamma_0(N)$ whose field of Fourier coefficients is $\Q(\sqrt2)$ ($N$ being a positive integer less than or equal to $500$), our table contains a curve of genus $2$ whose jacobian is conjecturably isogenous over $\Q$ to the abelian variety associated via the Shimura construction to $f$, so it is complete in this sense. Some of the curves are quotients of a modular curve by Atkin-Lehner involutions; in particular, their jacobians satisfy the generalised Shimura-Taniyama-Weil conjecture. 

We give two families of $2$-isogenies, each between jacobians of curves of genus $2$ with $\sqrt2$ multiplication. As we will see later, if $J$ is the jacobian of a curve of genus $2$, $A$ is an abelian variety, $\varepsilon$ is a $\sqrt2$ multiplication of $J$, and $\phi$ is an isogeny from $J$ to $A$ such that $\ker\phi$ is contained in $\ker\varepsilon$, then $A$ has a principal polarisation, and so is probably the jacobian of a curve of genus $2$; our families make this observation explicit. Also, as we will see later, the jacobian of a curve of genus $2$ whose endomorphisms are precisely the multiplication by $n$ endomorphisms (this is the generic case) is not $2$-isogenous to the jacobian of another curve of genus $2$, so the existence of the $\sqrt2$ multiplications is important here.

We give a family of curves of genus $2$ such that the ring of endomorphisms of the jacobian of a curve in it contains a maximal order in the division quaternion algebra $\lb(\frac{2,-3}{\Q}\rb)$; more precisely, it contains the endomorphisms $\varepsilon$ and $\eta$, where $\varepsilon$ is a $\sqrt2$ multiplication and $\eta$ is a $\frac{1+\sqrt{-3}}{2}$ multiplication such that $\eta\varepsilon+\varepsilon\eta=\varepsilon$. If $Y^2=F(X)$ is a curve in this family, then there is a $9$-isogeny defined over $\Q$ between the jacobians of the curves $Y^2=F(X),Y^2=-3F(X)$. As we will see later, there is an interesting connection with newforms; if $f$ is a newform such that the abelian variety associated to it via the Shimura construction is isogenous over $\Q$ to the jacobian of a curve in this family, then, letting $\sigma$ be the non-trivial automorphism of $\Q(\sqrt2)$, $\sigma(f)$ is the twist of $f$ by the Dirichlet character $\Z\rightarrow\C,n\mapsto\lb(\begin{array}{c}n\\3\end{array}\rb)$.

We give explicit examples of other isogenies between the jacobians of pairs of curves with $\sqrt2$ multiplication. Each isogeny is derived by regarding the jacobians as complex tori (via Abel-Jacobi isomorphisms), finding the appropriate morphism of complex tori, and then finding the appropriate morphism of abelian varieties.       

Throughout this paper, $K$ will be a subfield of $\C$, and $\overline{K}$ will be the subfield of $\C$ which is an algebraic closure of $K$. The field of definition of an object is defined to be the smallest field between $K$ and $\overline{K}$ over which it is defined. 

I would like to thank B. Birch for his supervision of my doctoral research, and various people for their support and insights, especially  J. Wilson, J. Merriman, A. Brumer, H. Cohn, and M. Muller, for his private communication (\cite{muller}). I would also like to thank the EPSRC for their financial support during my time at Kent.
\section{Homomorphisms between Jacobians of Curves}\label{isog}  
Throughout this paper, we will be concerned with homomorphisms between jacobians of curves of genus $2$. In this section, we recall standard facts about homomorphisms between jacobians of curves which will be useful later on, starting with curves of general genus (describing explicitly the Rosati involution w.r.t the canonical polarisation on the jacobian of a curve, called simply the Rosati involution in this paper, using the curve rather than the jacobian), and ending with curves of genus $2$ (associating to a homomorphism a $2$ by $2$ matrix, useful for handling compositions of homomorphisms). The material on curves of general genus is not simplified by restricting to curves of genus $2$, so it is described for curves of general genus even though it is only used for curves of genus $2$. It is taken from~\cite{LB}, Ch. 11, \S 5, although we use divisors instead of line bundles. 
\bigskip\\
Let $C_1,C_2$ be curves defined over $\overline{K}$. We define the following notation:
\begin{itemize}
\item $\pi_i$ is the projection from $C_1\times C_2$ to $C_i$.
\item $\tau$ is the involution from $C_1\times C_2$ to $C_2\times C_1$ defined by $\tau((P,Q))=(Q,P)$.  
\item $\lambda_{P}$ is the embedding from $C_2$ to $C_1\times C_2$ defined by $\lambda_P(Q)=(P,Q)$;
\item $P_i$ is a point on $C_i$;
\item $\Theta_i$ is the theta divisor $\{[P-P_i]\mid P\in C_i\}$ on the jacobian of $C_i$;
\item $\alpha_i$ is the embedding from $C_i$ to its jacobian defined by $\alpha_i(P)=[P-P_i]$;
\item $t_{i,x}$ is the translation by $x$ map on the jacobian of $C_i$;
\item $\Phi_{\Theta_i}$ is the map from the jacobian of $C_i$ to its dual defined by $\Phi_{\Theta_i}(x)=t_{i,x}^*([\Theta_i])-[\Theta_i]$ (the canonical polarisation).    
\item $\mu_i$ is the addition map from the square of the jacobian of $C_i$ to the jacobian of $C_i$;
\end{itemize}   
To any divisor class $[D]$ on $C_1\times C_2$, we can associate a homomorphism $\phi_{[D]}$ from the jacobian of $C_1$ to the jacobian of $C_2$ as follows:
\begin{equation} 
\phi_{[D]}\lb(\sum_{P\in C_1}n_P[P]\rb)=\sum_{P\in C_1}n_P\lambda_{P}^*([D]).\label{div}
\end{equation}
Moreover, any homomorphism arises in this way, and if $\phi_{[D]}$ is zero, then there are divisor classes $[D_1],[D_2]$ on $C_1,C_2$ respectively such that $[D]=\pi_1^*([D_1])+\pi_2^*([D_2])$. 
\bigskip\\  
We have the following proposition:
\begin{prop}\label{proprich}
$\Phi_{\Theta_1}\phi_{\tau^*([D])}\phi_{[D]}=\phi^*_{[D]}\Phi_{\Theta_2}\phi_{[D]}$. If $C_1=C_2$, then $\phi_{\tau^*([D])}=\phi^\dagger_{[D]}$, where $\dagger$ is the Rosati involution on the jacobian of $C_1$.
\end{prop}
\pf Firstly, we have
\begin{equation}
\phi_{[D]}=\phi_{(\phi_{[D]}\alpha_1\times\alpha_2)^*\mu_2^*([-\Theta_2])};\label{ident}
\end{equation}
indeed, for all points $P$ on $C_1$, we have 
\begin{eqnarray*}
\tabspace\phi_{(\phi_{[D]}\alpha_1\times\alpha_2)^*\mu_2^*([-\Theta_2])}\alpha_1(P)\\
\tabeq\phi_{(\phi_{[D]}\alpha_1\times\alpha_2)^*\mu_2^*([-\Theta_2])}([P-P_1])\\
\tabeq\lambda_P^*(\phi_{[D]}\alpha_1\times\alpha_2)^*\mu_2^*([-\Theta_2])-\lambda_{P_1}^*(\phi_{[D]}\alpha_1\times\alpha_2)^*\mu_2^*([-\Theta_2])\\
\tabeq\alpha_2^*t^*_{2,\phi_{[D]}\alpha_1(P)}([-\Theta_2])-\alpha_2^*t^*_{2,\phi_{[D]}\alpha_1(P_1)}([-\Theta_2])\\
\tabeq\alpha_2^*([\Theta_2]-t^*_{2,\phi_{[D]}\alpha_1(P)}([\Theta_2]))\\
\tabeq-\alpha_2^*\Phi_{\Theta_2}\phi_{[D]}\alpha_1(P)\\
\tabeq\phi_{[D]}\alpha_1(P),
\end{eqnarray*} 
the final equality following from~\cite{LB}, Ch. 11 Proposition 3.5. Secondly, we have
\begin{equation}
\phi_{\tau^*([D])}=\Phi^{-1}_{\Theta_1}\phi^*_{[D]}\Phi_{\Theta_2};\label{ident77}
\end{equation}
indeed, for all points $P$ on $C_2$, we have
\begin{eqnarray*}
\tabspace\phi_{\tau^*([D])}\alpha_2(P)\\
\tabeq\phi_{\tau^*(\phi_{[D]}\alpha_1\times\alpha_2)^*\mu_2^*([-\Theta_2])}\alpha_2(P)\\
\tabeq\phi_{(\alpha_2\times\phi_{[D]}\alpha_1)^*\mu_2^*([-\Theta_2])}\alpha_2(P)\\
\tabeq\phi_{(\alpha_2\times\phi_{[D]}\alpha_1)^*\mu_2^*([-\Theta_2])}(P-P_2)\\
\tabeq\lambda^*_P(\alpha_2\times\phi_{[D]}\alpha_1)^*\mu_2^*([-\Theta_2])-\lambda^*_{P_2}(\alpha_2\times\phi_{[D]}\alpha_1)^*\mu_2^*([-\Theta_2])\\
\tabeq\alpha^*_1\phi^*_{[D]}t^*_{2,\alpha_2(P)}([-\Theta_2])-\alpha^*_1\phi^*_{[D]}t^*_{2,\alpha_2(P_2)}([-\Theta_2])\\
\tabeq\alpha^*_1\phi^*_{[D]}([\Theta_2]-t^*_{2,\alpha_2(P)}([\Theta_2]))\\
\tabeq-\alpha^*_1\phi^*_{[D]}\Phi_{\Theta_2}\alpha_2(P)\\
\tabeq\Phi^{-1}_{\Theta_1}\phi^*_{[D]}\Phi_{\Theta_2}\alpha_2(P),
\end{eqnarray*}
the first equality following from (\ref{ident}) and the final equality following from~\cite{LB}, Ch. 11 Proposition 3.5.

The proposition now follows from (\ref{ident77}) (we recall that, if $C_1=C_2$, then $\phi^\dagger_{[D]}=\Phi_{\Theta_1}^{-1}\phi^*_{[D]}\Phi_{\Theta_1}$).
\finpf
\bigskip\\
>From now on in this section, we concentrate on curves of genus $2$, more explicitly.
\bigskip\\ 
Let $C_1,C_2$ be the curves
\begin{equation}
C_1\colon Y^2=F_1(X), C_2\colon T^2=F_2(Z),\label{c}
\end{equation}
where $F_1$ and $F_2$ are polynomials of degree $5$ or $6$ defined over $K$.

Let $P_1$ be a point on $C_1$ defined over $\overline{K}$. To any homomorphism $\phi$ from the jacobian of $C_1$ to the jacobian of $C_2$ defined over $K$, we can associate rational functions $R_1,\dots,R_4,S_1,\dots,S_4$ of $x$ defined over the field of definition of $P_1$ such that $\phi$ is given explicitly by
\[[(x,y)-P_1]\mapsto[(z_1,t_1)-(z_2,-t_2)],\] 
where $z_1,z_2$ are the zeros w.r.t. $x$ of
\begin{equation}
z^2+(R_3(x)+R_4(x)y)z+R_1(x)+R_2(x)y,\label{genxdef}
\end{equation} 
and $t_k$ equals  
\begin{equation} 
t_k=(S_3(x)+S_4(x)y)z_k+S_1(x)+S_2(x)y;\label{genydef}
\end{equation} 
this follows from considering the embedding of the jacobian of a curve of genus $2$ into $\p^{15}$ described in~\cite{CF}, Ch. 2.

Let $\iota_1$ be the hyperelliptic involution on $C_1$. If $\phi$ kills $[P_1-\iota_1(P_1)]$ (which clearly happens when $P_1$ is a Weierstrass point, but can happen when $P_1$ is not, as we will see later), then $R_2,R_4,S_1,S_3$ are all zero. 
\bigskip\\
Finally in this section, we will associate a $2$ by $2$ matrix to a homomorphism between jacobians of curves of genus $2$, called its {\em differential matrix}; this association will be used frequently in this paper. To do this, we will discuss jacobians of curves of genus $2$ in terms of complex tori.

Firstly, we will recall the Abel-Jacobi isomorphism from the jacobian of a curve of genus $2$ to a complex torus; this will allow us to discuss the former in terms of the latter. Let $C$ be a curve of genus $2$ defined over $\C$, of the form $Y^2=F(X)$, where $F$ is a polynomial of degree $5$ or $6$ defined over $\C$ (a basis for its vector space over $\C$ of holomorphic differentials is $\frac{dX}{Y},\frac{X\:dX}{Y}$), and define a lattice $\Lambda$ of $\C^2$ as follows:
\[\Lambda=\lb\{\int_c\lb(\frac{dx}{y},\frac{x\:dx}{y}\rb)\mid c\in H_1(C,\Z)\rb\}.\]  
It is well-known that the map 
\[\Jac(C)\rightarrow\C^2/\Lambda,[P-P']\mapsto\int_{P'}^{P}\lb(\frac{dx}{y},\frac{x\:dx}{y}\rb)\]
is an isomorphism, called the {\em Abel-Jacobi isomorphism} (see~\cite{LB}, Ch. 11 \S 1). 

Before discussing homomorphisms between jacobians of curves of genus $2$ in terms of complex tori, we will define some complex-analytic notation which will be useful in this section and in later sections, stating standard and well-known results (see~\cite{LB}, Ch. 1). Let $c_1,\dots,c_4$ be a symplectic basis of $H_1(C,\Z)$. The matrix
\[\Pi=\lb(\begin{array}{cccc}\int_{c_1}\frac{dx}{y}&\int_{c_2}\frac{dx}{y}&\int_{c_3}\frac{dx}{y}&\int_{c_4}\frac{dx}{y}\\\int_{c_1}\frac{x\:dx}{y}&\int_{c_2}\frac{x\:dx}{y}&\int_{c_3}\frac{x\:dx}{y}&\int_{c_4}\frac{x\:dx}{y}\end{array}\rb)\]
is called a {\em period matrix} of the jacobian of $C$. For use in \S\S~\ref{quat} and~\ref{other}, we define a positive definite hermitian form on $\C^2$ as follows: 
\[H(\underline{z},\underline{w})=2i\underline{z}^t\lb(\overline{\Pi}\lb(\begin{array}{cccc}0&0&-1&0\\0&0&0&-1\\1&0&0&0\\0&1&0&0\end{array}\rb)\Pi^t\rb)^{-1}\overline{\underline{w}}.\]
$H$ is the hermitian form on $\C^2$ corresponding to the canonical polarisation. It is independent of the choice of the symplectic basis.
   
We will now discuss homomorphisms between jacobians of curves of genus $2$ in terms of complex tori, stating standard and well-known results (see~\cite{LB}, Ch. 1). Let $C_1,C_2$ be curves of the form (\ref{c}), let $\Lambda_i$ be the lattice corresponding to $C_i$, and let $\Pi_i$ be a period matrix of the jacobian of $C_i$; then, a homomorphism $\phi$ from the jacobian of $C_1$ to the jacobian of $C_2$ induces via the Abel-Jacobi isomorphisms a homomorphism from $\C^2/\Lambda_1$ to $\C^2/\Lambda_2$, to which can be associated a $2$ by $2$ matrix $A_\phi$ with complex elements (independent of the choice of the symplectic bases) and a $4$ by $4$ matrix $R_\phi$ (dependent on the choice of the symplectic bases) with integer elements such that  
\begin{equation}
A_\phi\Pi_1=\Pi_2R_\phi\label{pereqn} 
\end{equation} 
(the notation $A_\phi,R_\phi$ comes from the complex-analytic terminology `analytic representation' and `rational representation'). Moreover, $A_\phi$ is uniquely determined by $\phi$; via this association, composition of homomorphisms is compatible with composition of matrices, and the homomorphism is zero if and only if the matrix is zero.

We are now in a position to associate a $2$ by $2$ matrix to a homomorphism between jacobians of curves of genus $2$, called its {\em differential matrix}, as described in the following proposition:
\begin{prop}\label{dmprop}
Let $C_1,C_2$ be curves of the form (\ref{c}), and let $\phi$ be a homomorphism from the jacobian of $C_1$ to the jacobian of $C_2$ given explicitly by (\ref{genxdef}) and (\ref{genydef}). Then,
\begin{equation}
\lb(\begin{array}{c}\sum_{k=1}^2\frac{dz_k}{t_k}\\\sum_{k=1}^2\frac{z_k\:dz_k}{t_k}\end{array}\rb)=A_\phi\lb(\begin{array}{c}\frac{dx}{y}\\\frac{x\:dx}{y}\end{array}\rb),\label{diffmat}
\end{equation}
where $A_\phi$ is as defined above; we call $A_\phi$ the {\em differential matrix} of $\phi$.

Via this association, the actions of the absolute Galois group of $K$ on homomorphisms defined over $\overline{K}$ and on matrices defined over $\overline{K}$ are compatible, composition of homomorphisms is compatible with composition of matrices, and the homomorphism is zero if and only if the matrix is zero (in particular, this association induces an explicit embedding from the ring of endomorphisms defined over $K$ of the jacobian of a curve of genus $2$ defined over $K$ to the ring of $2$ by $2$ matrices defined over $K$).  

The differential matrix can be computed in steps as follows: 
\begin{eqnarray*}
G_1\tabeq R_2'+R_2\frac{F_1'}{2F_1},G_2=R_4'+R_4\frac{F_1'}{2F_1},\\\\
H_1\tabeq2(R_1R_3'+G_2R_2F_1)-R_1'R_3-G_1R_4F_1,H_2=2(R_2R_3'+G_2R_1)-R_1'R_4-G_1R_3,\\
H_3\tabeq R_3'R_3+G_2R_4F_1-2R_1',H_4=R_3'R_4+G_2R_3-2G_1,\\
H_5\tabeq R_3^2+R_4^2F_1-4R_1,H_6=2R_3R_4-4R_2,\\\\
I_1\tabeq\frac{H_1H_5-H_2H_6F_1}{H_5^2-H_6^2F_1},I_2=\frac{H_2H_5-H_1H_6}{H_5^2-H_6^2F_1},I_3=\frac{H_3H_5-H_4H_6F_1}{H_5^2-H_6^2F_1},I_4=\frac{H_4H_5-H_3H_6}{H_5^2-H_6^2F_1},\\\\
J_1\tabeq-R_1I_3-R_2I_4F_1,J_2=-R_1I_4-R_2I_3,\\
J_3\tabeq I_1-R_3I_3-R_4I_4F_1,J_4=I_2-R_3I_4-R_4I_3,\\\\
L_3\tabeq S_3^2+S_4^2F_1,L_4=2S_3S_4,L_5=S_1S_3+S_2S_4F_1,L_6=S_1S_4+S_2S_3,\\
L_1\tabeq S_1^2+S_2^2F_1+R_1L_3+R_2L_4F_1-R_3L_5-R_4L_6F_1,\\
L_2\tabeq2S_1S_2+R_1L_4+R_2L_3-R_3L_6-R_4L_5,\\\\
M_3\tabeq I_1S_3+I_3S_1+(I_2S_4+I_4S_2)F_1,M_4=I_1S_4+I_2S_3+I_3S_2+I_4S_1,\\
M_5\tabeq2(I_3S_3+I_4S_4F_1),M_6=2(I_3S_4+I_4S_3),\\
M_1\tabeq2(I_1S_1+I_2S_2F_1)+R_1M_5+R_2M_6F_1-R_3M_3-R_4M_4F_1,\\
M_2\tabeq2(I_1S_2+I_2S_1)+R_1M_6+R_2M_5-R_3M_4-R_4M_3,\\\\
N_3\tabeq J_1S_3+J_3S_1+(J_2S_4+J_4S_2)F_1,N_4=J_1S_4+J_2S_3+J_3S_2+J_4S_1,\\
N_5\tabeq2(J_3S_3+J_4S_4F_1),N_6=2(J_3S_4+J_4S_3),\\
N_1\tabeq2(J_1S_1+J_2S_2F_1)+R_1N_5+R_2N_6F_1-R_3N_3-R_4N_4F_1,\\
N_2\tabeq2(J_1S_2+J_2S_1)+R_1N_6+R_2N_5-R_3N_4-R_4N_3,
\end{eqnarray*}
\begin{eqnarray}
\frac{(A_\phi)_{11}+(A_\phi)_{12}x}{y}\tabeq\sum_{k=1}^2\frac{1}{t_k}\frac{dz_k}{dx}=\frac{(L_1M_2-L_2M_1)F_1+(L_1M_1-L_2M_2F_1)y}{(L_1^2-L_2^2F_1)y},\label{ide1}\\
\frac{(A_\phi)_{21}+(A_\phi)_{22}x}{y}\tabeq\sum_{k=1}^2\frac{z_k}{t_k}\frac{dz_k}{dx}=\frac{(L_1N_2-L_2N_1)F_1+(L_1N_1-L_2N_2F_1)y}{(L_1^2-L_2^2F_1)y}\nonumber
\end{eqnarray}
(for brevity, we denote $F(x)$ by $F$, where $F$ is a rational function, and differentiation w.r.t. $x$ by $'$).
\end{prop}
\pf First, we will show that (\ref{diffmat}) holds for some $2$ by $2$ matrix $A_\phi$ with complex elements; later we will show that $A_\phi$ is the differential matrix of $\phi$. From (\ref{genxdef}) and (\ref{genydef}), we have
\begin{equation}
\lb(\begin{array}{c}\sum_{k=1}^2\frac{dz_k}{t_k}\\\sum_{k=1}^2\frac{z_k\:dz_k}{t_k}\end{array}\rb)=\lb(\begin{array}{c}\frac{(A_1(x)+A_2(x)y)dx}{y}\\\frac{(A_3(x)+A_4(x)y)dx}{y}\end{array}\rb)\label{diff2}
\end{equation}
for some rational functions $A_1,A_2,A_3,A_4$ (divide both elements of both sides by $dx$ and use the fact that the LHS remains invariant after swapping $(z_1,t_1)$ and $(z_2,t_2)$); to show that (\ref{diffmat}) holds for some $2$ by $2$ matrix $A_\phi$, we need to show that both differentials in the RHS of (\ref{diff2}) are holomorphic, or, equivalently, that $A_1,A_3$ are linear polynomials and that $A_2,A_4$ are zero. Let $\mcC$ be a path on $C_1$, and let $\mcC_1,\mcC_2$ be the paths on $C_2$ induced by (\ref{genxdef}) and (\ref{genydef}). Then, we have
\begin{equation}
\lb(\begin{array}{c}\sum_{k=1}^2\int_{\mcC_k}\frac{dz_k}{t_k}\\\sum_{k=1}^2\int_{\mcC_k}\frac{z_k\:dz_k}{t_k}\end{array}\rb)=\lb(\begin{array}{c}\int_\mcC\frac{(A_1(x)+A_2(x)y)dx}{y}\\\int_\mcC\frac{(A_3(x)+A_4(x)y)dx}{y}\end{array}\rb).\label{in}
\end{equation}
Since $\frac{dz}{t},\frac{z\:dz}{t}$ are holomorphic differentials, all four integrals in the LHS of (\ref{in}) are not infinity, so both integrals in the RHS of (\ref{in}) are not infinity. This holds for all paths $\mcC$ on $C_1$, so both differentials in the RHS of (\ref{diff2}) are holomorphic, as required.

Now, we will show that $A_\phi$ is the differential matrix of $\phi$. This follows from (\ref{in}); indeed, we see that there is a $4$ by $4$ matrix $R_\phi$ with integer elements satisfying (\ref{pereqn}), since $\phi$ induces an element of $H_1(C_2,\Z)$ from an element of $H_1(C_1,\Z)$, via (\ref{genxdef}) and (\ref{genydef}).
\bigskip\\
Consider the second paragraph of the theorem. We remarked when discussing homomorphisms between complex tori that, via this association, composition of homomorphisms is compatible with composition of matrices, and the homomorphism is zero if and only if the matrix is zero. From (\ref{diffmat}), we see that the actions of the absolute Galois group of $K$ on homomorphisms defined over $\overline{K}$ and on matrices defined over $\overline{K}$ are compatible as well.
\bigskip\\
Consider the final paragraph of the theorem. This follows from evaluating both elements of the LHS of (\ref{diffmat}) divided by $dx$ step by step, making the denominators as simple as possible. More precisely, we have
\begin{eqnarray*}
\frac{dz}{dx}\tabeq\frac{R_1+G_1y+z(R_3+G_2y)}{-R_3-R_4y-2z}=\frac{H_1+H_2y+z(H_3+H_4y)}{H_5+H_6y}=I_1+I_2y+z(I_3+I_4y),\\
\frac{z\:dz}{dx}\tabeq J_1+J_2y+z(J_3+J_4y),\\
\prod_{k=1}^2t_k\tabeq L_1+L_2y,\sum_{k=1}^2t_{3-k}\frac{dz_k}{dx}=M_1+M_2y,\sum_{k=1}^2z_kt_{3-k}\frac{dz_k}{dx}=N_1+N_2y,
\end{eqnarray*}
\begin{eqnarray*}
\sum_{k=1}^2\frac{1}{t_k}\frac{dz_k}{dx}\tabeq\frac{\sum_{k=1}^2t_{3-k}\frac{dz_k}{dx}}{\prod_{k=1}^2t_k}=\frac{M_1+M_2y}{L_1+L_2y}=\frac{(L_1M_2-L_2M_1)F_1+(L_1M_1-L_2M_2F_1)y}{(L_1^2-L_2^2F_1)y},\\
\sum_{k=1}^2\frac{z_k}{t_k}\frac{dz_k}{dx}\tabeq\frac{\sum_{k=1}^2z_kt_{3-k}\frac{dz_k}{dx}}{\prod_{k=1}^2t_k}=\frac{N_1+N_2y}{L_1+L_2y}=\frac{(L_1N_2-L_2N_1)F_1+(L_1N_1-L_2N_2F_1)y}{(L_1^2-L_2^2F_1)y}.
\end{eqnarray*}
We note that
\[\frac{(L_1M_2-L_2M_1)F_1}{L_1^2-L_2^2F_1},\frac{(L_1N_1-L_2N_2)F_1}{L_1^2-L_2^2F_1}\]
are linear polynomials, and that
\[L_1M_1-L_2M_2F_1,L_1N_1-L_2N_2F_1\]
are zero, from (\ref{ide1}). We have not attempted to prove this directly.
\finpf 
\bigskip\\
To compute the differential matrix, it suffices to compute it at two points $(x,y)$, giving two linear equations for its top row elements, and two linear equations for its bottom row elements. This observation considerably eases its computation.        

As will be demonstrated frequently in this paper, the explicit  embedding from the ring of endomorphisms defined over $K$ of the jacobian of a curve of genus $2$ defined over $K$ to the ring of $2$ by $2$ matrices defined over $K$ allows us add and compose two endomorphisms by adding and composing their differential matrices, which is much easier.      
\section{The Richelot dual}\label{rich} 
A key tool in the derivation of our family of curves of genus $2$ whose jacobians have $\sqrt2$ multiplication is the Richelot dual of a curve of genus $2$ and isogenies associated to it, which we now discuss.  
\bigskip\\
Let $C$ be the curve
\[Y^2=\Delta\prod_{i=0}^2G_i(X),\]
defined over $K$, where $\Delta\in K^*$, the $G_i$'s are quadratics defined over $\overline{K}$ individually and over $K$ as a set, and the sextic $\Delta\prod_{i=0}^2G_i(X)$ has no multiple zeros. Let $g_{ij}\in\overline{K}$ such that $G_i(X)=\sum_{j=0}^2g_{ij}X^j$.
\bigskip\\
Suppose for now that $\det(g_{uv})$ is non-zero (later on, we will look at the case where $\det(g_{uv})$ is zero), and define the {\em Richelot dual of $C$} to be the curve
\[\det(g_{uv})T^2=\Delta\prod_{i=0}^2H_i(Z),\]
defined over $K$, where $H_i(Z)=G'_{i+1}(Z)G_{i+2}(Z)-G'_{i+2}(Z)G_{i+1}(Z)$ (here, and in similar places throughout this paper, we will take addition to mean addition modulo $3$). 
\bigskip\\
We will now define homomorphisms defined over $\overline{K}$ from the jacobian of $C$ to the jacobian of the Richelot dual of $C$, and from the jacobian of the Richelot dual of $C$ to the jacobian of $C$, in the manner of \S~\ref{isog}. We will show that they are isogenies defined over $K$ using Proposition~\ref{dmprop}. 

There is a homomorphism $\rho$ defined over $\overline{K}$ from the jacobian of $C$ to the jacobian of the Richelot dual of $C$, called the {\em Richelot isogeny}, given explicitly by
\[[(x,y)-P_0]\mapsto[(z_1,t_1)-(z_2,-t_2)],\]
where $P_0$ is a Weierstrass point on $C$ whose $x$-coordinate is a zero of $G_0$, $z_1,z_2$ are the zeros w.r.t. $z$ of
\begin{equation}
G_1(x)H_1(z)+G_2(x)H_2(z),\label{xdef}
\end{equation}
and $t_k$ satisfies
\begin{equation}
yt_k=\Delta G_1(x)H_1(z_k)(x-z_k).\label{ydef}
\end{equation}
By Lemma 2.6.1 in~\cite{bending}, $(z_1,t_1)$ and $(z_2,t_2)$ are points on the Richelot dual of $C$.

Also, there is a homomorphism $\varrho$ defined over $\overline{K}$ from the jacobian of the Richelot dual of $C$ to the jacobian of $C$, called the {\em Richelot dual isogeny}, given explicitly by
\[[(z,t)-Q_0]\mapsto[(x_1,y_1)-(x_2,-y_2)],\]
where $Q_0$ is a Weierstrass point on the Richelot dual of $C$ whose $z$-coordinate is a zero of $H_0$, $x_1,x_2$ are the zeros w.r.t. $x$ of (\ref{xdef}), and $y_k$ satisfies
\begin{equation}
ty_k=\Delta G_1(x_k)H_1(z)(x_k-z).\label{ydef2}
\end{equation} 
By Lemma 2.6.1 in~\cite{bending}, $(x_1,y_1)$ and $(x_2,y_2)$ are points on $C$.

As we will explain below, the Richelot isogeny and the Richelot dual isogeny are isogenies defined over $K$, and we have
\begin{equation}
\varrho\rho=[2],\rho\varrho=[2],\label{dualrel}
\end{equation}
justifying the terminology `Richelot dual isogeny'.
\bigskip\\
We will now state some standard facts about the Richelot isogeny and the Richelot dual isogeny, which will be used in this paper.

Firstly, via some computation, we see that the differential matrices of the Richelot isogeny (resp. the Richelot dual isogeny) are the identity (resp. twice the identity). So they are isogenies defined over $K$, and we have (\ref{dualrel}), by Proposition~\ref{dmprop}.

Secondly, define $P_{ij},Q_{ij}$ ($i\in\{0,1,2\},j\in\{0,1\}$) to be the Weierstrass points on $C$ (resp. its Richelot dual) whose $x$-coordinates (resp. $z$-coordinates) are the zeros of $G_i$ (resp. $H_i$). Then, we have
\begin{eqnarray} 
\ker\rho=\{0,[P_{00}-P_{01}],[P_{10}-P_{11}],[P_{20}-P_{21}]\},\label{kerrich}\\
\ker\varrho=\{0,[Q_{00}-Q_{01}],[Q_{10}-Q_{11}],[Q_{20}-Q_{21}]\}, 
\nonumber
\end{eqnarray} 
and
\begin{eqnarray}
\rho([P_{00}-P_{10}])=[Q_{20}-Q_{21}],\rho([P_{00}-P_{20}])=[Q_{10}-Q_{11}],\label{imrich}\\
\varrho([Q_{00}-Q_{10}])=[P_{20}-P_{21}],\varrho([Q_{00}-Q_{20}])=[P_{10}-P_{11}],\nonumber
\end{eqnarray}   
as we now explain. Via some computation, we see that (\ref{kerrich}) with $=$ replaced by $\supseteq$ holds, and that (\ref{imrich}) holds. We see that the degrees of the Richelot isogeny and the Richelot dual isogeny are at least $4$. However, from (\ref{dualrel}), they are both $4$, establishing (\ref{kerrich}). 
\bigskip\\
The following proposition will be vital when looking at curves of genus $2$ with $\sqrt2$ multiplication:
\begin{prop}\label{pol} 
The pullback of the canonical polarisation on the jacobian of $C$ by the Richelot dual isogeny is twice the canonical polarisation on the jacobian of the Richelot dual of $C$.
\end{prop}  
\pf If $D$ is the divisor on the product of the Richelot dual of $C$ with $C$ induced by (\ref{xdef}) and (\ref{ydef2}), then the Richelot dual isogeny is $\phi_{[D]}$ and the Richelot isogeny is $\phi_{\tau^*([D])}$ (adopting notation from \S~\ref{isog}); this follows from (\ref{div}) and the definition of $\tau$, bearing in mind that (\ref{ydef}) is (\ref{ydef2}) with $x_k,z,y_k,t$ replaced by $x,z_k,y,t_k$ respectively. The proposition now follows from Proposition~\ref{proprich} and (\ref{dualrel}).
\finpf  
\bigskip\\
The Richelot dual, the Richelot isogeny and the Richelot dual isogeny are related to an algorithm derived by Richelot in the nineteenth century to calculate, with great rapidity and precision, integrals of the form
\[\int_u^{u'}\frac{S(x)dx}{\sqrt{-P(x)}},\int_v^{v'}\frac{S(x)dx}{\sqrt{-P(x)}},\int_w^{w'}\frac{S(x)dx}{\sqrt{-P(x)}},\]
where $u,u',v,v',w,w'$ are strictly increasing real numbers, $P$ is the monic sextic whose zeros are $u,u',v,v',w,w'$, and $S$ is a polynomial of degree $\leq1$. For further details, see~\cite{bostmestre}, \S A2.
%
\bigskip\\ 
We now look at the case where $\det(g_{uv})$ is zero. We have the following lemma:
\begin{lemma}\label{split}
Suppose that $\det(g_{uv})$ is zero. Then, the jacobian of $C$ has two abelian subvarieties defined over $\overline{K}$ which are elliptic curves; call them $E_+,E_-$, and there is an isogeny $\rho$ defined over $\overline{K}$ from the jacobian of $C$ to $E_+\times E_-$ such that
\begin{equation}
\ker\rho=\{0,[P_{00}-P_{01}],[P_{10}-P_{11}],[P_{20}-P_{21}]\}.\label{ellker}
\end{equation}
Moreover, there is an isogeny $\varrho$ defined over $\overline{K}$ from $E_+\times E_-$ to the jacobian of $C$ such that
\begin{equation}
\varrho\rho=[2],\label{dualrel2}
\end{equation}
and the pullback of the canonical polarisation on the jacobian of $C$ by $\varrho$ is twice the canonical polarisation on $E_+\times E_-$. 
\end{lemma}
\pf Since $\det(g_{uv})$ is zero, there is an automorphism $\alpha$ of $C$, swapping for each $i$ the Weierstrass points corresponding to $G_i$, such that $\sigma^{-1}\alpha\sigma=\pm\alpha$ if $\alpha$ is an automorphism of $\overline{K}$ over $K$. 

Define $E_+,E_-$ to be $\im([1]+\alpha),\im([1]-\alpha)$ respectively (here, $\alpha$ denotes the induced automorphism of the jacobian of $C$); we need to show that they are elliptic curves, or, equivalently, that the endomorphisms $[1]+\alpha,[1]-\alpha$ are neither zero nor surjective. They are not zero, since $[1]\neq\alpha,-\alpha$, and they are not surjective, since either both or none are surjective, being defined over $K$ as a pair, and $([1]+\alpha)([1]-\alpha)=[0]$.

We define $\rho$ as follows:
\[\rho(P)=(([1]+\alpha)(P),([1]-\alpha)(P)).\]
Now the kernel of $\rho$ consists precisely of the points on the jacobian of $C$ which are sent by $[1],\iota,-\iota$ to the same point. Since $\iota$ is an automorphism, all of these points are killed by the multiplication by $2$ endomorphism. Since $\alpha$ swaps for each $i$ the Weierstrass points corresponding to $G_i$, we have (\ref{ellker}), as required.   

We define $\varrho$ as follows:
\[\varrho((P_+,P_-))=P_++P_-.\]
Then, we have
\[\varrho\rho(P)=\varrho(([1]+\alpha)(P)+([1]-\alpha)(P))=[2]P,\]
establishing (\ref{dualrel2}).

Finally, since $\alpha$ is fixed by the Rosati involution w.r.t. the canonical polarisation (being induced from an automorphism of $C$), $[1]+\alpha,[1]-\alpha$ are primitive, and $([1]+\alpha)^2=[2]([1]+\alpha),([1]-\alpha)^2=[2]([1]-\alpha)$, the pullback of the canonical polarisation on the jacobian of $C$ by the inclusion homomorphism from $E_+,E_-$ to it is twice the canonical polarisation on $E_+,E_-$ respectively, by~\cite{LB}, Ch. 5 Criterion 3.4. So the pullback of the canonical polarisation on the jacobian of $C$ by $\varrho$ is twice the canonical polarisation on $E_+\times E_-$, by~\cite{LB}, Ch. 12 Lemma 1.6.
\finpf 
\section{Curves with $\protect\sqrt2$ Multiplication}\label{root2} 
In this section, we give a family of curves of genus $2$ defined over $K$ whose
jacobians have $\sqrt2$ multiplication defined over $K$ and fixed by the Rosati
involution; moreover, any curve with these properties is isomorphic over $K$ to a curve in our family. More precisely, we prove the following theorem:
\begin{theorem}\label{Thm3.1.1}
Let $C$ be a curve of genus $2$ defined over $K$, and suppose that the jacobian of $C$ has a $\sqrt2$ multiplication $\varepsilon$ defined over $K$ and fixed by the Rosati involution. 

Then, there is an isomorphism $\psi$ defined over $K$ from $C$ to a curve of the form
\begin{equation}Y^2=\Delta\prod_{i=0}^2G_i(X), (G_i(X)=X^2-\alpha_i
X+P\alpha_i^2+Q\alpha_i+R),\label{root2curve}
\end{equation}
where
\begin{eqnarray}
\prod_{i=0}^2(X-\alpha_i)=X^3+AX^2+BX+C,\label{cubic4}\\
R=4P,B=\frac{Q(PA-Q)+4P^2+1}{P^2},C=\frac{4(PA-Q)}{P},\label{param}
\end{eqnarray}
and
\begin{equation}
\Delta,P,Q,A\in K,\Delta\neq0\neq P.\label{conditions}
\end{equation}
Moreover, $\psi\varepsilon\psi^{-1}=\pm\iota^{-1}\rho$, where $\rho$ is the Richelot isogeny, and $\iota$ is the isomorphism from $C$ to its Richelot dual defined by
\begin{equation}
(x,y)\mapsto\lb(\frac{2}{x},\frac{4y}{x^3}\rb)\label{root2isom}  
\end{equation}
(here, $\psi$ and $\iota$ also denote the induced isomorphisms between the jacobians).    

Conversely, the jacobian of a curve $C$ of the
form (\ref{root2curve}),
where (\ref{cubic4}), (\ref{param}) and (\ref{conditions}) hold, has a
$\sqrt2$ multiplication $\varepsilon$ defined over $K$ and fixed by the Rosati involution, namely $\iota^{-1}\rho$.  
\end{theorem}
The assumption that $\varepsilon$ is fixed by the Rosati involution is weak, since it is satisfied if the ring of endomorphisms defined over $K$ is generated by $\varepsilon$ (which is usually the case), and the Rosati involution cannot negate a $\sqrt2$ multiplication (see~\cite{LB}, Ch. 5 Theorem 1.8). In particular, this is the case if $K=\Q$ and the jacobian of $C$ satisfies the generalised Shimura-Taniyama-Weil conjecture, since then its ring of endomorphisms defined over $\Q$ embeds into the ring of $2$ by $2$ matrices defined over $\Q$, by Proposition~\ref{dmprop}, but is not equal to it, a consequence of satisfying the conjecture. A family of principally polarised abelian surfaces defined over $\C$ given as complex tori, whose members are generically isomorphic to jacobians of curves of genus $2$ and have infinitely many $\sqrt2$ multiplications but none fixed by the Rosati involution, is given in~\cite{bending}, \S 3.8.      

\S~\ref{table} contains a table of curves defined over $\Q$ whose jacobians have $\sqrt2$ multiplication defined over $\Q$ and are simple over $\Q$, and so satisfy the conditions of the generalised Shimura-Taniyama-Weil conjecture; all the curves are isomorphic over $\Q$ to a curve from Theorem~\ref{Thm3.1.1}, as expected.
\bigskip\\
\pf {\bf(of Theorem~\ref{Thm3.1.1})} The theorem consists of two parts; the second part is the `Conversely' statement, and the first part is the rest. We will prove the first before the second, as the first is the harder.
\bigskip\\
We will now prove the first part of the theorem. To start with, we note that, by the standard procedure based on the Riemann-Roch theorem, $C$ is isomorphic over $K$ to a curve of the form $Y^2=F(X)$, where $F$ is a polynomial of degree $5$ or $6$ defined over $K$ (see~\cite{CF}, Ch. 1). So we may assume, without loss of generality, that $C$ is a curve of this form.
\bigskip\\
The material in \S~\ref{rich} will be crucial; bearing (\ref{kerrich}), (\ref{imrich}) and (\ref{ellker}) in mind, it will be useful to discover how $\varepsilon$ behaves on the $2$-torsion on the jacobian of $C$. In fact, we will show that there is a labelling $P_{ij},i\in\{0,1,2\},j\in\{0,1\}$ of the Weierstrass points on $C$ such that
\begin{equation}
\ker\varepsilon=\{0,[P_{00}-P_{01}],[P_{10}-P_{11}],[P_{20}-P_{21}]\}\label{ker} 
\end{equation} 
and 
\begin{equation}
\varepsilon([P_{00}-P_{10}])=[P_{20}-P_{21}],\varepsilon([P_{00}-P_{20}])=[P_{10}-P_{11}].\label{2tors}
\end{equation} 

To do this, we will use Weil pairings (see~\cite{mil}, \S 16 for the definitions). If $m$ is a positive integer, we define
\[e_m\colon\Jac(C)[m]\times\Jac(C)[m]\rightarrow\mu_m\]
to be the Weil pairing corresponding to $m$ times the canonical polarisation on the jacobian of $C$. Now, since $\varepsilon$ is fixed by the Rosati involution, we have
\begin{equation}
e_2(x,y)=e_1(\varepsilon(x),\varepsilon(y))=1\:\forall x,y\in\ker\varepsilon\label{ident3}
\end{equation}
(the pullback of the canonical polarisation by $\varepsilon$ is twice the canonical polarisation), and
\begin{equation}
e_2(x,\varepsilon(x))=1\:\forall x,y\in\Jac(C)[2].\label{ident13}
\end{equation}   
It is well-known that $e_2([P_a-P_b],[P_c-P_d])$ equals $1$ if and only if the sets $\{a,b\},\{c,d\}$ are equal or have an empty intersection (see~\cite{CF}, Ch. 3 \S 3). So, from (\ref{ident3}), we deduce that (\ref{ker}) holds, and, from (\ref{ident13}), we deduce that (\ref{2tors}) holds.
\bigskip\\
Let $G_i(X)$ be the monic quadratic, defined over $\overline{K}$, whose zeros are the $x$-coordinates of $P_{i0}$ and $P_{i1}$, and let $\rho$ (resp. $\varrho$) be the Richelot isogeny (resp. the Richelot dual isogeny). Adopting notation from \S~\ref{rich}, the crux of the proof of the first part of the theorem is to show that the Richelot dual of $C$ exists and that there is an isomorphism defined over $K$ from $C$ to it which sends  $P_{ij}$ to $Q_{ij}$ (after, for each $i$, swapping $Q_{i0},Q_{i1}$ if necessary); then, the first part of the theorem will follow from computations performed in~\cite{bending}, Ch. 3. 
\bigskip\\
Suppose, for a contradiction, that the Richelot dual of $C$ does not exist. Then, $\det(g_{uv})$ is zero. Taking notation from Lemma~\ref{split}, there is an isomorphism $\iota$ defined over $\overline{K}$ from the jacobian of $C$ to $E_+\times E_-$ such that
\[\varepsilon=\iota^{-1}\rho=\varrho\iota;\]
this follows from (\ref{ker}), (\ref{ellker}) and (\ref{dualrel2}). Let $\Theta$ be a theta divisor on the jacobian of $C$, and let $D$ be the divisor on the jacobian of $C$ which is the pullback of the canonical divisor on $E_+\times E_-$ by $\iota$. By Lemma~\ref{split} and the above equality and the fact that $\varepsilon$ is fixed by the Rosati involution, the pullback of the canonical polarisation on $E_+\times E_-$ by $\iota$ is the canonical polarisation on the jacobian of $C$, so the polarisations corresponding to $\Theta$ and $D$ are equal. In particular, $\Theta\cdot D=\Theta^2$; we will obtain our contradiction by showing that this equation does not hold. 

Now $D$ contains the trivial point and six non-trivial points of order $2$ (three on each elliptic curve); in particular, it contains $0,[P_0-P_1],[P_0-P_2]$ for some Weierstrass points $P_0,P_1,P_2$ on $C$. Replacing $\Theta$ by a translate if necessary, we may assume that $\Theta$ contains $0,[P_0-P_1],[P_0-P_2]$ as well, which implies that $\Theta\cdot D\geq3$. But it is a well-known consequence of the Riemann-Roch theorem that $\Theta^2=2$ (see~\cite{mil}, Theorem 13.3), giving our contradiction.          
\bigskip\\
We have now established that the Richelot dual of $C$ exists. Letting $\rho$ (resp. $\varrho$) be the Richelot isogeny (resp. the Richelot dual isogeny), there is an isomorphism $\iota$ defined over $K$ from the jacobian of $C$ to the jacobian of its Richelot dual such that
\begin{eqnarray}
\varepsilon\tabeq\iota^{-1}\rho,\label{eqn20}\\
\tabeq\varrho\iota;\label{eqn21}
\end{eqnarray}
this follows from the fact that $\varepsilon$ and $\rho$ are defined over $K$, (\ref{ker}), (\ref{kerrich}) and (\ref{dualrel}). By Lemma~\ref{pol} and the second equality and the fact that $\varepsilon$ is fixed by the Rosati involution, $\iota$ respects the canonical polarisations. We are now in a position to show that there is an isomorphism defined over $K$ from $C$ to its Richelot dual which sends $P_{ij}$ to $Q_{ij}$ (after, for each $i$, swapping $Q_{i0},Q_{i1}$ if necessary).
\bigskip\\
Applying Torelli's theorem to $\iota$ (see~\cite{mil2}, \S\S 12 and 13 for the theorem and its proof), there is an isomorphism $\iota$ defined over $K$ from $C$ to its Richelot dual, such that 
\begin{equation}
\iota([P-P'])=[\iota(P)-\iota(P')]\label{eqn10}
\end{equation}
for all points $P,P'$ on $C$ (justifying the notation $\iota$ for both isomorphisms).    

Considering the image under (\ref{eqn20}) of $[P_{10}-P_{20}]$ and using (\ref{2tors}) and (\ref{imrich}), we see that $\iota([P_{00}-P_{01}])=[Q_{00}-Q_{01}]$, and so, from (\ref{eqn10}), that
\begin{equation}
\{\iota(P_{00}),\iota(P_{01})\}=\{Q_{00},Q_{01}\}.\label{eqn11} 
\end{equation} 
Considering the image under (\ref{eqn21}) of $[P_{00}-P_{10}]$ and using (\ref{2tors}) and (\ref{imrich}), we see that $\iota([P_{00}-P_{01}])$ is one of the four points $[Q_{0j_1}-Q_{1j_2}]$, and so, from (\ref{eqn10}) and (\ref{eqn11}), that
\begin{equation}
\iota(P_{1j})=Q_{1j},\label{eqn12}
\end{equation}
swapping $Q_{10},Q_{11}$ if necessary. Similarly,
\begin{equation}
\iota(P_{2j})=Q_{2j},\label{eqn13}
\end{equation}
swapping $Q_{20},Q_{21}$ if necessary. From (\ref{eqn11}), (\ref{eqn12}) and (\ref{eqn13}), we see, swapping $Q_{00},Q_{01}$ if necessary, that $\iota$ sends $P_{ij}$ to $Q_{ij}$, as required. We have now proved the first part of the theorem.  
\bigskip\\
We will now prove the second part of the theorem. The differential matrices of $\rho$ (resp. $\iota$) are the identity (resp. $\lb(\begin{array}{cc}0&-\frac{1}{2}\\-1&0\end{array}\rb)$), so the differential matrix of $\iota^{-1}\rho$ is $\lb(\begin{array}{cc}0&-1\\-2&0\end{array}\rb)$. We deduce that $\iota^{-1}\rho$ is a $\sqrt2$ multiplication defined over $K$, by Proposition~\ref{dmprop}. It is fixed by the Rosati involution by Proposition~\ref{pol} and the fact that the pullback of a theta divisor on the jacobian of $C$ by $\iota^{-1}$ is a theta divisor on the jacobian of its Richelot dual. We have now proved the second part of the theorem.    
\finpf 
\bigskip\\
The following proposition can be used to show that the jacobians of some (but not all) of the curves in the table contained in \S~\ref{table} are simple over $\Q$, and will be used when we look at curves of genus $2$ whose jacobians have $\lb(\frac{2,-3}{\Q}\rb)$ multiplication in \S~\ref{quat}:
\begin{prop}\label{split2}
Let $C$ be a curve of genus $2$ defined over $K$. Suppose that its jacobian has a $\sqrt2$ multiplication $\varepsilon$ defined over $K$, and that the kernel of $\varepsilon$ is the only subgroup of order $4$ of the $2$-torsion of its jacobian defined over $K$.

Then, its jacobian is simple over $K$.
\end{prop}
\pf See~\cite{bending}, Proposition 3.3.1. The proof essentially uses the fact that the Galois group of the $2$-torsion of the square of an elliptic curve defined over $K$ is relatively small, whereas the Galois group of the $2$-torsion of an abelian surface defined over $K$ is relatively big if there is only one subgroup of order $4$ defined over $K$.   
\finpf
\section{$2$-Isogenies between Curves with $\protect\sqrt2$ Multiplication}  
In this section, we give two families of $2$-isogenies defined over $K$, each between jacobians of curves of genus $2$ defined over $K$ with $\sqrt2$ multiplication defined over $K$ and fixed by the Rosati involution. More precisely, we prove the following theorems:
\begin{theorem}\label{isogthm} 
Suppose that $\Delta,U,V\in K^*,W\in K$, and define 
\[F_1(X)=\Delta(X^2+UX+V)\prod_{i=1}^2\lb(X^2-\alpha_i X-\frac{U\alpha_i+4}{V}\rb),\]
where
\[\prod_{i=1}^2(X-\alpha_i)=X^2+\frac{W(U-V)(U+V)+4(V^2+4)}{4U}X+W.\]
We suppose throughout this theorem that $F_1$ has six distinct zeros, and define $C_1$ to be the curve $C_1\colon Y^2=F_1(X)$.

Also, we suppose throughout this theorem that $V\neq2$, and define
\[F_2(Z)=\Delta(Z^2+U'Z+V')\prod_{i=1}^2\lb(Z^2-\alpha'_i Z-\frac{U'\alpha_i+4}{V'}\rb),\]
where
\[\prod_{i=1}^2(Z-\alpha'_i)=Z^2+\frac{W'(U'-V')(U'+V')+4({V'}^2+4)}{4U'}Z+W',\]
and
\begin{eqnarray*}U'\tabeq\frac{2(V-U+2)}{U-2},\\
V'\tabeq\frac{2(U-V)}{U-2},\\
W'\tabeq\frac{W(U-V)(V^3+V^2(U-8)+4(V(U+1)-U))}{(U-V)^2(W(V-1)-4V)+4(V-2)^2}\\
\tabplus\frac{4(V^4-8(V^3-V^2(U+1)+V(U^2-2(U-2))-2))}{(U-V)^2(W(V-1)-4V)+4(V-2)^2},\\
\Delta'\tabeq\frac{\Delta(U-2)((U-V)^2(W(V-1)-4V)+4(V-2)^2)}{4UV(V-2)}\\
\end{eqnarray*}
(the denominators of $U',V',W'$ disappear when $F_2$ is expanded, and so can be zero).

Then, $F_2$ has six distinct zeros (possibly including $\infty$); we define $C_2$ to be the curve $C_2\colon Y^2=F_2(Z)$.
\bigskip\\
The jacobians of $C_1$ and $C_2$ both have a $\sqrt2$ multiplication defined over $K$ and fixed by the Rosati involution, killing the points $D_1$ (resp. $D_2$) of order $2$ corresponding to $X^2+UX+V$ (resp. $(U-2)Z^2+2(V-U+2)Z+2(U-V)$), and the other two points of order $2$ corresponding to factors of $F_1$ (resp. $F_2$).

Moreover, there is a $2$-isogeny $\pi_1$ defined over $K$ from the jacobian of $C_1$ to the jacobian of $C_2$ killing $D_1$, and a $2$-isogeny $\pi_2$ defined over $K$ from the jacobian of $C_2$ to the jacobian of $C_1$ killing $D_2$; they are given explicitly by
\[[(x,y)-P_1]\mapsto[(z_1,t_1)-(z_2,-t_2)],[(z,t)-P_2]\mapsto[(x_1,t_1)-(x_2,-t_2)],\]
where $P_1$ (resp. $P_2$) is a Weierstrass point on $C_1$ (resp. $C_2$) supporting $D_1$ (resp. $D_2$), $z_1,z_2$ (resp. $x_1,x_2$) are the zeros w.r.t. $z$ (resp. $x$) of
\[\sum_{i=0}^2\sum_{j=0}^2\Phi_{ij}x^iz^j,\]
and $t_k,y_k$ equal
\begin{eqnarray*}
t_k\tabeq y\frac{(V-2)(\Psi_1(x)z_k+\Psi_2(x))}{(x^2+Ux+V)\lb(\sum_{i=0}^2\Phi_{i2}x^i\rb)^2},\\
y_k\tabeq t\frac{(V-2)(\Psi_3(z)x_k+\Psi_4(z))}{((U-2)z^2+2(V-U+2)z+2(U-V))\lb(\sum_{j=0}^2\Phi_{2j}z^j\rb)^2},
\end{eqnarray*}
with $\Phi_{ij}$ being the element in the $i^{th}$ row and the $j^{th}$ column of the $3$ by $3$ matrix
\[\lb(\begin{array}{ccc}2V(U-2)&2V(4-U)&2(U-2(V-1))\\2(U(U-2)+V(V-2))&2(U(V-U+2)-(V-2)^2)&U(U-V)-2(V-2)\\V(U+V-4)&2(2V-U)&U-V\end{array}\rb),\]
and $\Psi_1(X),\Psi_2(X),\Psi_3(Z),\Psi_4(Z)$ being defined as follows:
\begin{eqnarray*}
\Psi_1(X)\tabeq2(U^2-V(V+2))X^3\\
\tabminus(UV^2+(U+4)(6V-U(U+2)))X^2\\
\tabminus(V^2(4V-U^2-16)+(U+2)(2V(U+10)+U^3-4(2U^2-U+2)))X\\
\tabminus2(4V^3-(U+2)((U+V^2)(U+2)-V(U^2-U+8))),\\
\Psi_2(X)\tabeq(8V-U(U+V)(V-U+2))X^3\\
\tabplus((U+V)(U-V)(U^2-2V)-4(3V^2-8V+U^2))X^2\\
\tabplus2(4V^2(V-3)-(U+2)(V^2U+V(U-8)-U(U-1)(U-2)))X\\
\tabplus2V(4V^2-(U+2)(V(U+2)-U^2+2(U-2))).\\ 
\Psi_3(Z)\tabeq(8V-U(U+V)(V-U+2))Z^3\\
\tabminus(2V^3+(U-2)(U+2)V^2+2(U-4)^2V-U^2(U^2-4(U-2)))Z^2\\
\tabminus(3UV^3-2(U-4)V^2-(U^2(U+10)-32(U-1))V+2U^2(U-2)^2)Z\\
\tabminus2(V^3(V-2)+U(U-2)(4V-U^2)),\\
\Psi_4(Z)\tabeq4V(2-V)Z^3\\
\tabminus2((3U-8)V^2+2(U^2-U+4)V-U^2(U+2))Z^2\\
\tabminus2V(2V^2+(U^2-4(U-3))V+U^2(U-6))Z\\
\tabminus2V(UV^2+2(U-4)V-U^2(U-2)).
\end{eqnarray*}
\end{theorem} 
\begin{theorem}\label{isogthm2}
Suppose that $\Delta\in K^*,W\in K$, and define
\[F_1(X)=\Delta(X^2-2)\prod_{i=1}^2(X^2-\alpha_i X+2),\]
where 
\[\prod_{i=1}^2(X-\alpha_i)=X^2+WX+8.\]
We suppose throughout this theorem that $F_1$ has six distinct zeros, and define $C_1$ to be the curve $C_1\colon Y^2=F_1(X)$.

Also, we define
\[F_2(Z)=\Delta(Z^2-2)\prod_{i=1}^2(Z^2-\alpha'_i Z+2),\]
where 
\[\prod_{i=1}^2(Z-\alpha'_i)=Z^2+W'Z+8,\]
and 
\[W'=\frac{16(6-W)}{16-3W},\Delta'=\Delta(16-3W)\]
(the denominator of $W'$ disappears when $F_2$ is expanded, and so can be zero).

Then $F_2$ has six distinct zeros (possibly including $\infty$); we define $C_2$ to be the curve $C_2\colon Y^2=F_2(Z)$.
\bigskip\\
The jacobians of $C_1$ and $C_2$ both have a $\sqrt2$ multiplication defined over $K$ and fixed by the Rosati involution, killing the points $D_1$ (resp. $D_2$) of order $2$ corresponding to $X^2+UX+V$ (resp. $(U-2)Z^2+2(V-U+2)Z+2(U-V)$), and the other two points of order $2$ corresponding to factors of $F_1$ (resp. $F_2$).

Moreover, there is a $2$-isogeny $\pi_1$ defined over $K$ from the jacobian of $C_1$ to the jacobian of $C_2$ killing $D_1$, and a $2$-isogeny $\pi_2$ defined over $K$ from the jacobian of $C_2$ to the jacobian of $C_1$ killing $D_2$; they are given explicitly by
\[[(x,y)-(\sqrt2,0)]\mapsto[(z_1,t_1)-(z_2,-t_2)],[(z,t)-(\sqrt2,0)]\mapsto[(x_1,t_1)-(x_2,-t_2)],\]
where $z_1,z_2$ (resp. $x_1,x_2$) are the zeros w.r.t. $z$ (resp. $x$) of
\[3x^2z^2+4(x^2z+xz^2)+2(x^2+z^2)+16xz+8(x+z)+12,\]
and $t_k,y_k$ equal
\begin{eqnarray*}
t_k\tabeq32y\frac{(x+1)(x^2-4x-2)z_k-x^3-2x^2-6x-4}{(x^2-2)(3x^2+4x+2)^2},\\
y_k\tabeq4t\frac{z(5z^2+12z+6)x_k+2(2z^3+7z^2+8z+2)}{(z^2-2)(3z^2+4z+2)^2}.
\end{eqnarray*}
\end{theorem}
It can happen that one value of $\Delta$ and three values of $(U,V,W)$ give rise to the same curve $C_1$ in Theorem~\ref{isogthm}. This is because the $\sqrt2$ multiplication of its jacobian kills three points of order $2$ defined over $K$; a particular value of $(U,V,W)$ corresponds to a particular point. The three curves $C_2$ are in general distinct.       
\bigskip\\
Suppose that $C$ is a curve of genus $2$ defined over $K$ whose jacobian has a $\sqrt2$ multiplication defined over $K$ which kills a point of order $2$ defined over $K$. Then, $C$ is isomorphic over $K$ to a curve $C_1$ from Theorem~\ref{isogthm} (possibly with $V=2$) or from Theorem~\ref{isogthm2}. This follows from some computation, as we now explain. By Theorem~\ref{Thm3.1.1}, $C$ is isomorphic over $K$ to a curve of the form
\begin{equation}
Y^2=\Delta(X^2+UX+V)\prod_{i=0}^2(X^2-\alpha_iX+Q\alpha_i+R)\label{iscurve}
\end{equation}
for some $U,V,Q,R\in K$ and quadratic conjugates $\alpha_1,\alpha_2\in\overline{K}$; if 
\[G_0(X)=X^2+UX+V,G_i(X)=X^2-\alpha_iX+Q\alpha_i+R,\]
then there is an isomorphism from (\ref{iscurve}) to its Richelot dual defined by (\ref{root2isom}). There are two cases to consider: the case $U\neq0$ leads to a curve $C_1$ from Theorem~\ref{isogthm} (possibly with $V=2$), and the case $U=0$ leads to a curve $C_1$ from Theorem~\ref{isogthm2}. 
\bigskip\\
The jacobian of a curve of genus $2$ defined over $K$ whose endomorphisms defined over $K$ are precisely the multiplication by $n$ endomorphisms (this is the generic case) is not $2$-isogenous over $K$ to the jacobian of another curve of genus $2$ defined over $K$, as we now explain. Let $\phi$ be a $2$-isogeny from the jacobian $J_1$ to the jacobian $J_2$ defined over $K$, and let $\Phi_i$ be the canonical polarisation on $J_i$. Then, there is an endomorphism of degree $4$ defined over $K$ of $J_1$, namely $\Phi^{-1}_1\phi^*\Phi_2\phi$. However, the multiplication by $n$ endomorphism of $J_1$ has degree $n^4$, and so is not $4$. 

However, if $C$ is a curve of genus $2$ defined over $K$ whose jacobian has a $\sqrt2$ multiplication $\varepsilon$ defined over $K$ and fixed by the Rosati involution, $A$ is an abelian surface defined over $K$, and $\phi$ is a $2$-isogeny defined over $K$ from the jacobian of $C$ to $A$ such that $\ker\phi$ is contained in $\ker\varepsilon$, then, as we now explain, $A$ has a principal polarisation defined over $K$ (generically coming from a curve of genus $2$ on $A$) and a $\sqrt2$ multiplication defined over $K$ and fixed by its Rosati involution; our theorems make this observation explicit. To do this, we will use Weil pairings, as we did for part of our proof of Theorem~\ref{Thm3.1.1}. Let $\Phi_J$ be the canonical polarisation on the jacobian of $C$, and define 
\[e\colon\Jac(C)[[2]+\varepsilon]\times\Jac(C)[[2]+\varepsilon]\rightarrow\mu_2\]
to be the Weil pairing corresponding to $\Phi_J([2]+\varepsilon)$, which is a polarisation since $[2]+\varepsilon$ is fixed by the Rosati involution and the zeros of the characteristic polynomial of its differential matrix are positive (see~\cite{LB}, Ch. 5 Theorem 2.4). Since the Weil pairing is bilinear and alternating, $e$ is trivial on $\ker\phi$, so there is a principal polarisation $\Phi_A$ on $A$ such that $\Phi_J([2]+\varepsilon)=\phi^*\Phi_A\phi$ (see~\cite{mil}, Proposition 16.8). Moreover, if $\hat{\phi}$ is the $8$-isogeny defined over $K$ from $A$ to the jacobian of $C$ such that $[2]=\hat{\phi}\phi$, then $\hat{\phi}\varepsilon\phi$ is twice a $\sqrt2$ multiplication of $A$ defined over $K$ (since $\ker\phi$ is contained in $\ker\varepsilon$) and fixed by the Rosati involution w.r.t. $\Phi_A$ (via a bit of algebraic manipulation). 
\bigskip\\
When $V=2$, the jacobian of a curve $C_1$ from Theorem~\ref{isogthm} is $4$-isogenous over at worst a quadratic extension of $K$ to the square of an elliptic curve, as we explain below, so its quotient abelian variety by the subgroup generated by $P$ could well be the direct sum of two abelian subvarieties which are elliptic curves defined over at worst a quadratic extension of $K$ rather than the jacobian of a curve of genus $2$ defined over $K$, which would justify the condition $V\neq2$ in the theorem. If $\alpha\in\overline{K}$ is satisfies
\begin{equation}
(U-2)\alpha^2+4\alpha+U+2=0,\label{alpharel}
\end{equation}
and $E$ is the elliptic curve
\begin{eqnarray}
E\colon T^2\tabeq\frac{\Delta(\alpha+1)}{U}(Z+1)\label{ellcur}\\
\tabtimes[-32(4(U-3)\alpha+U^2-2U-4)(Z^2-6Z+1)\nonumber\\
\tabplus(U-2)^2((U^2-12)\alpha-2(U+2))W(Z^2+1)\nonumber\\
\tabplus2((U^4-4U^3-8U^2-16U+144)\alpha-2(U^3+6U^2-20U-24))WZ]\nonumber
\end{eqnarray}
(whose coefficients of $1,Z^3$ and of $Z,Z^2$ are equal), then there is a $4$-isogeny defined over $\Q(\alpha)$ from the jacobian of $C_1$ to the square of $E$, given explicitly by
\[[(x,y)-(x',y')]\mapsto\lb((z^2,t)-({z'}^2,t'),\lb(\frac{1}{z^2},\frac{t}{z^3}\rb)-\lb(\frac{1}{{z'}^2},\frac{t'}{{z'}^3}\rb)\rb),\]
where
\[z=\frac{x-\alpha-1}{\alpha x+\alpha-1},t=\frac{32y(U^3-8U^2+4U+32+\alpha(U-4)(U^2+4U-20))}{((U-2)(\alpha x+\alpha-1))^3},\]
$z',t'$ having analogous definitions. 
\bigskip\\
\pf {\bf(of Theorem~\ref{isogthm})} Most of our proof involves straightforward but tedious computation. We will briefly discuss two points: the existence, kernels and the fields of definition of the $\sqrt2$ multiplications, and the kernels of the isogenies, which require some mathematics.

Consider the $\sqrt2$ multiplications. The curves $C_1,C_2$ can be treated analogously, so we will just consider $C_1$. After some computation, we see that there is an isomorphism $\iota$ defined over $K$ from $C_1$ to its Richelot dual defined by (\ref{root2isom}). Letting $\iota$ denote the induced isomorphism between the jacobians, we seee that $\iota^-1\rho$ is a $\sqrt2$ multiplication $\varepsilon_1$ defined over $K$ of the jacobian of $C_1$ with the required kernel (cf. the final paragraph of our proof of Theorem~\ref{Thm3.1.1}).

Consider the isogenies. After some computation, we see that their differential matrices are
\[A_{\pi_1}=\lb(\begin{array}{cc}-1&0\\0&-1\end{array}\rb),A_{\pi_2}=\lb(\begin{array}{cc}-2&1\\2&-2\end{array}\rb),\]
and so, by Proposition~\ref{dmprop}, that $\pi_2\pi_1$ is $[2]\pm\varepsilon_1$, an isogeny of degree $4$ (by multiplying the differential matrices above). We deduce that they have the required kernels.
\finpf
\bigskip\\
\pf {\bf(of Theorem~\ref{isogthm2})} Our proof of Theorem~\ref{isogthm2} is analogous to our proof of Theorem~\ref{isogthm}, the computations less tedious. The differential matrices of the isogenies are
\[A_{\pi_1}=\lb(\begin{array}{cc}-\frac{1}{2}&-\frac{1}{2}\\1&\frac{1}{2}\end{array}\rb),A_{\pi_2}=\lb(\begin{array}{cc}0&2\\-4&0\end{array}\rb).\]
\finpf
\section{Curves with $\lb(\frac{2,-3}{\Q}\rb)$ Multiplication}\label{quat} 
In this section, we give a family of curves of genus $2$ defined over $K$ such that the ring of endomorphisms defined over $\overline{K}$ of the jacobian of a curve in it contains a maximal order in the division quaternion algebra $\lb(\frac{2,-3}{\Q}\rb)$. More precisely, we prove the following theorem:
\begin{theorem}\label{quatthm}
Let $\Delta,N\in K^*$ such that $N^2\neq-1$. Then, the ring of endomorphisms defined over $\overline{K}$ of the jacobian of the curve
\begin{eqnarray*}
C\tabcolon\\
Y^2\tabeq\Delta((N-1)X^6-6NX^5+3(N+1)X^4-8N^2X^3+3(N-1)X^2+6NX+N+1)
\end{eqnarray*}
contains the endomorphisms $\varepsilon$ and $\eta$, where $\varepsilon$ is a $\sqrt2$ multiplication and $\eta$ is a
$\frac{1+\sqrt{-3}}{2}$ multiplication such that $\eta\varepsilon+\varepsilon\eta=\varepsilon$ (the ring generated by $\varepsilon$ and $\eta$ is a maximal order in the division quaternion algebra $\lb(\frac{2,-3}{\Q}\rb)$). Moreover, $\varepsilon$ is defined over $K$, $\eta$ is defined over $K(\sqrt{-3})$, $\sigma^{-1}\eta\sigma+\eta=[1]$ if $K$ does not contain $\sqrt{-3}$ and $\sigma$ is the non-trivial automorphism of $K(\sqrt{-3})$ over $K$, and, under the Rosati involution, $\varepsilon$ maps to itself and $\eta$ maps
to $[2]+\varepsilon-\eta([3]+[2]\varepsilon)$.

If $N\in\Q$ and $N\neq\pm1$, then the ring of endomorphisms defined over $\overline{K}$ is generated by the endomorphisms $\varepsilon$ and $\eta$.
\end{theorem}
A family equivalent to ours was derived by K. Hashimoto and N. Murabayashi, their curves being defined over quadratic fields (see~\cite{HM}, Theorem 1.3). 

Since the ring of endomorphisms defined over $\Q$ of the jacobian of a curve of genus $2$ defined over $\Q$ embeds into the ring of $2$ by $2$ matrices defined over $\Q$, by Proposition~\ref{dmprop}, an order in a division quaternion algebra cannot embed into the former. The ring of endomorphisms defined over $\Qbar$ of the jacobian of a curve in our family with $N\in\Q$ is generated by one defined over $\Q$ and one defined over $\Q(\sqrt{-3})$.

There is an interesting connection with newforms. Let $f$ be a newform such that the abelian variety associated to it via the Shimura contruction is isogenous over $\Q$ to the jacobian of a curve in our family with $N\in\Q$. By our theorem, the jacobians of the curves
\begin{eqnarray*}
Y^2\tabeq\Delta((N-1)X^6-6NX^5+3(N+1)X^4-8N^2X^3+3(N-1)X^2+6NX+N+1),\\
Y^2\tabeq-3\\
\tabtimes\Delta((N-1)X^6-6NX^5+3(N+1)X^4-8N^2X^3+3(N-1)X^2+6NX+N+1)
\end{eqnarray*}
are isogenous over $\Q$, so the twist of $f$ by the Dirichlet character $\Z\rightarrow\C,n\mapsto\lb(\begin{array}{c}n\\3\end{array}\rb)$ is equal to either $f$ or $\sigma(f)$, where $\sigma$ is the non-trivial automorphism of $\Q(\sqrt2)$. In fact, it is the latter; if it was the former, then the jacobian of our curve would be isogenous over $\Qbar$ to the square of an elliptic curve with complex multiplication by $\Q(\sqrt{-3})$ (see~\cite{ribet}, \S 5), which is impossible (if $N\neq\pm1$, then, by our theorem, the jacobian of our curve is not isogenous over $\Qbar$ to the square of an elliptic curve with complex multiplication, and if $N=\pm1$, then the jacobian of our curve is isogenous over $\Qbar$ to the square of an elliptic curve with complex multiplication by $\Q(\sqrt{-6})$) (see~\cite{bending}, \S 5.1).
\bigskip\\              
\pf {\bf(of Theorem~\ref{quatthm})} The first paragraph of the theorem can be proved by explicit computation, as follows. Without loss of generality, we may assume that $N\neq1$ and that $\Delta=N-1$. We will define the endomorphism $\varepsilon$ and an endomorphism $\psi$, and then the endomorphism $\eta$ in terms of them.  
\bigskip\\
Consider the endomorphism $\varepsilon$. For $i=0,1,2$, define $G_i$ as
\begin{equation}
G_i(X)=X^2(\alpha_i-1)-2\alpha_i^2X+\alpha_i+1,\label{expandcurve} 
\end{equation}
where the $\alpha_i$'s are the cube roots of $N$; then, the $x$-coordinates of the Weierstrass points on $C$ are the zeros of the $G_i$'s. There is an isomorphism $\iota$ from $C$ to its Richelot dual defined by
\[(x,y)\mapsto\lb(\frac{1-x}{x+1},\frac{4y}{(x+1)^3}\rb),\]
which induces an isomorphism $\iota$ between the jacobians. We define $\varepsilon$ to be $\iota^{-1}\rho$, where $\rho$ is the Richelot isogeny; we see that its differential matrix is 
\begin{equation}
A_\varepsilon=\lb(\begin{array}{cc}-1&-1\\-1&1\end{array}\rb),\label{root2dm}
\end{equation}
(via Proposition~\ref{dmprop}), and so, by Proposition~\ref{dmprop}, that it is a $\sqrt2$ multiplication defined over $K$. It is also fixed by the Rosati involution (cf. the final paragraph of our proof of Theorem~\ref{Thm3.1.1}).
\bigskip\\
Consider the endomorphism $\psi$. We define $\psi$ to be given explicitly by
\[[(x,y)-\infty^+]\mapsto[(z_1,t_1)-(z_2,-t_2)],\]
where $z_1,z_2$ are the zeros of
\begin{eqnarray*}
\tabspace z^2+(\sqrt{-3}-1)(2x+1+\sqrt{-3})^{-1}x\\
\tabplus z(1+\sqrt{-3})(2x+1+\sqrt{-3})^{-1}(2x^2+2Nx-1+\sqrt{-3}(2Nx+1))^{-1}\\
\tabtimes(2(N+1)x^3-(2N-1)x^2-(6N+1)x-2N+\sqrt{-3}((4N+1)x^2+x)-2y),
\end{eqnarray*}
and $t_k$ equals
\begin{eqnarray*}
\tabspace(2x+1+\sqrt{-3})^{-2}(2x^2+2Nx-1+\sqrt{-3}(2Nx+1))^{-2}\\
\tabtimes[2(\sqrt{-3}-1)(2x^2+2Nx-1+\sqrt{-3}(2Nx+1))\\
\tabtimes\{(N-1)x(2x^4+(3N+2)x^3-3(3N+1)x^2-(7N+4)x-N-1\\
\tabplus\sqrt{-3}((3N+2)x^3+3(N+1)x^2-3Nx-N-1))\\
\tabplus y(2x^2+(N+1)x-1+\sqrt{-3}((N+1)x+1))\}\\
\tabplus4(1+\sqrt{-3})z_k\\
\tabtimes\{(N-1)(2x^7+2(N^2+5N+1)x^6-(18N^2+21N+10)x^5+(6N^2-46N-13)x^4\\
\tabplus(8N^3+48N^2+8N+3)x^3+2(21N^2+13N+4)x^2+3(N+1)(6N+1)x\\
\tabplus2N^2+6N+1\\
\tabplus\sqrt{-3}(2(N+1)^2x^6+(6N^2+9N+4)x^5-3(10N^2+6N+1)x^4\\
\tabplus(8N^3-12N^2-34N-9)x^3+2(3N^2-7N-1)x^2+3(2N^2-N+1)x+2N^2+1))\\
\tabplus y(2x^4-2(N^2+N-1)x^3+(N-4)(4N+1)x^2+(12N^2-7N-4)x+2N^2+4N-1\\
\tabplus\sqrt{-3}(-2(N^2-2N-1)x^3-(4N^2+3N-4)x^2-9Nx+2N^2-2N-1))\}];
\end{eqnarray*}
we see that its differential matrix is
\begin{equation}
A_\psi=\lb(\begin{array}{cc}0&\frac{1}{2}(-1+\sqrt{-3})\\\frac{1}{2}(-1-\sqrt{-3})&1\end{array}\rb)\label{root5dm}
\end{equation}
(via some computation), and so, by Proposition~\ref{dmprop}, that it is a $\frac{1+\sqrt{5}}{2}$ multiplication defined over $K(\sqrt{-3})$. $\psi$ was derived using the methods in \S~\ref{other}; we decided on a $\frac{1+\sqrt{5}}{2}$ multiplication for convenience, $5$ being the smallest possible discriminant of a real quadratic field.  
\bigskip\\
We will now define $\eta$ in terms of $\varepsilon$ and $\psi$, using complex-analytic theory.

Regard the jacobian of $C$ as a complex torus $T=\C^2/\Lambda$, let $H$ be the corresponding hermitian form on $\C^2$ as defined in \S~\ref{isog}, and let $M$ be the positive definite hermitian matrix such that $H(\underline{z},\underline{w})=\underline{z}^tM\overline{\underline{w}}$ for all $\underline{z},\underline{w}\in\C^2$. If $\phi$ is an endomorphism, $\phi^\dagger$ is its image under the Rosati involution, and $A_\phi,A_{\phi^\dagger}$ are their differential matrices, then it is a standard result that
\[H_i(A_\phi\underline{z},\underline{w})=H_i(\underline{z},A_{\phi^\dagger}\underline{w})\]
for all $\underline{z},\underline{w}\in\C^2$ (see~\cite{LB}, Ch. 5 Proposition 1.1), which implies that
\begin{equation}
A_{\phi^\dagger}=\overline{M^{-1}A^t_\phi M}.\label{adjrel}
\end{equation}

Since $\varepsilon$ is fixed by the Rosati involution, there are real numbers $s,t$ such that 
\[M=\lb(\begin{array}{cc}s-t&-t\\-t&s+t\end{array}\rb),\]
by (\ref{root2dm}) and (\ref{adjrel}) with $\phi=\varepsilon$. We will now show that there are coprime integers $s_0,t_0$ such that $(s:t)=(s_0:t_0)$, and that $s_0^2-2t_0^2$ divides $2$; then, we will show that $s_0^2-2t_0^2=1$, by showing firstly that $s_0^2-2t_0^2$ is positive and secondly that $s_0^2-2t_0^2\neq2$. Finally, we will establish the first paragraph of the theorem. 
\bigskip\\
Now the sum of $\psi$ and its image under the Rosati involution is the endomorphism  
\begin{equation}
\frac{1}{s^2-2t^2}[-2(s+t)t[1]-(s+2t)t\varepsilon+2s^2\psi+2st\psi\varepsilon]\label{endo}
\end{equation}
(by Proposition~\ref{dmprop}, (\ref{root2dm}), (\ref{root5dm}) and (\ref{adjrel}), the differential matrices of these endomorphisms are the same). Also, the ring generated by $\varepsilon$ and $\psi$ is a maximal order in the division quaternion algebra $\lb(\frac{2,-3}{\Q}\rb)$ by Proposition~\ref{dmprop} (consider their differential matrices), and the ring of endomorphisms defined over $\overline{K}$ of the jacobian of $C$ is a free abelian group of finite rank (see~\cite{mil}, Proposition 12.5), so the four coefficients in the endomorphism (\ref{endo}) are all integers. So there are coprime integers $s_0,t_0$ such that $(s:t)=(s_0:t_0)$, and that $s_0^2-2t_0^2$ divides $2$. 
\bigskip\\
Now $\varepsilon$ and (\ref{endo}) are fixed by the Rosati involution, so the endomorphism
\[\varepsilon+s_0([1]+\varepsilon-[2]\psi)+t_0([2]+\varepsilon-[2]\psi\varepsilon)\]
is. Its square is $[2+3(s_0^2-2t_0^2)]$, so $s_0^2-2t_0^2$ is positive (see~\cite{LB}, Ch. 5 Theorem 1.8).
\bigskip\\
Suppose that $s_0^2-2t_0^2=2$; we will show that this is in fact impossible. The idea is to show that there are two $\sqrt2$ multiplications with the same kernel and fixed by the Rosati involution, not the negatives of each other, and then that their ratio is induced by an automorphism of $C$, by Theorem~\ref{Thm3.1.1}, which will lead to the impossibility. 

Now the endomorphism
\[\frac{s_0+2t_0}{2}[1]+\frac{1+s_0+t_0}{2}\varepsilon-s_0\psi-t_0\psi\varepsilon\]
(the coefficients are integers since $s_0^2-2t_0^2=2$ and $s_0,t_0$ are coprime) is a $\sqrt2$ multiplication fixed by the Rosati involution, by Proposition~\ref{dmprop}, (\ref{root2dm}), (\ref{root5dm}) and (\ref{adjrel}) (consider its differential matrix) whose kernel is that of $\varepsilon$ (the coefficients of $\varepsilon$ and $\psi\varepsilon$ are integers and the coefficients of $[1]$ and $\psi$ are even integers since $s_0^2-2t_0^2=2$ and $s_0,t_0$ are coprime). So this endomorphism equals $\zeta\varepsilon$ for some automorphism $\zeta$. Since $\varepsilon$ and $\zeta\varepsilon$ are fixed by the Rosati involution, $\zeta$ is of the form
\begin{eqnarray*}
\tabspace[(x,y)-(x',y')]\\
\tabmapsto\lb[\lb(\frac{M_{11}x+M_{12}}{M_{21}x+M_{22}},\frac{My}{(M_{21}x+M_{22})^3}\rb)-\lb(\frac{M_{11}x'+M_{12}}{M_{21}x'+M_{22}},\frac{My'}{(M_{21}x'+M_{22})^3}\rb)\rb]
\end{eqnarray*}
for some $M_{11},M_{12},M_{21},M_{22},M\in\overline{K}$, by Theorem~\ref{Thm3.1.1}. 

Now, by Proposition~\ref{dmprop}, the differential matrix of $\zeta$ is
\[\lb(\begin{array}{cc}\frac{\sqrt{-3}}{4}s_0+\frac{1}{2}&\frac{\sqrt{-3}}{4}(-s_0-2t_0)\\\frac{\sqrt{-3}}{4}(-s_0+2t_0)&-\frac{\sqrt{-3}}{4}s_0+\frac{1}{2}\end{array}\rb),\] 
so, multiplying $M_{11},M_{12},M_{21},M_{22}$ by a non-zero constant, we have
\[M_{11}=\frac{\sqrt{-3}}{4}s_0-\frac{1}{2},M_{12}=\frac{\sqrt{-3}}{4}(s_0-2t_0),M_{21}=\frac{\sqrt{-3}}{4}(s_0+2t_0),M_{22}=-\frac{\sqrt{-3}}{4}s_0-\frac{1}{2}.\]
After some computation, using the fact that $s_0^2-2t_0^2=2$, we see that
\[U=V^{-1}WV,\] 
where $U,V,W$ are the M\"obius transformations 
\[x\mapsto\frac{M_{11}x+M_{12}}{M_{21}x+M_{22}},x\mapsto\frac{(s_0+2)x+s_0-2t_0}{(2-s_0)x+2t_0-s_0},x\mapsto\omega x\] 
respectively, with $\omega=-\frac{1}{2}(1+\sqrt{-3})$, a primitive third root of unity.

Let $x_{i0},x_{i1}$ be the zeros of $G_i$ in $\overline{K}$ ($G_i(X)$ being as in (\ref{expandcurve})). Since $\varepsilon$ and $\zeta\varepsilon$ are $\sqrt2$ multiplications, $U$ permutes the pairs $\{x_{00},x_{01}\},\{x_{10},x_{11}\},\{x_{20},x_{21}\}$. So there are elements $\kappa_1,\kappa_2$ of $\overline{K}$ such that $V$ sends $x_{00}$ to $\kappa_1$, $x_{01}$ to $\kappa_2$, $\{x_{10},x_{11}\}$ to $\{\omega\kappa_1,\omega\kappa_2\}$, and $\{x_{20},x_{21}\}$ to $\{\omega^2\kappa_1,\omega^2\kappa_2\}$. Now 
\[\lb|\begin{array}{ccc}1&\kappa_1+\kappa_2&\kappa_1\kappa_2\\1&\omega\kappa_1+\omega^2\kappa_2&\omega\kappa_1\omega^2\kappa_2\\1&\omega^2\kappa_1+\omega\kappa_2&\omega^2\kappa_1\omega\kappa_2\end{array}\rb|=0,\]
so 
\[\lb|\begin{array}{ccc}1&\kappa_1|_{V^{-1}}+\kappa_2|_{V^{-1}}&\kappa_1|_{V^{-1}}\kappa_2|_{V^{-1}}\\1&(\omega\kappa_1)|_{V^{-1}}+(\omega^2\kappa_2)|_{V^{-1}}&(\omega\kappa_1)|_{V^{-1}}(\omega^2\kappa_2)|_{V^{-1}}\\1&(\omega^2\kappa_1)|_{V^{-1}}+(\omega\kappa_2)|_{V^{-1}}&(\omega^2\kappa_1)|_{V^{-1}}(\omega\kappa_2)|_{V^{-1}}\end{array}\rb|=0,\]
(here, $|_{V^{-1}}$ means the image under $V^{-1}$), so
\begin{equation}
\lb|\begin{array}{ccc}1&x_{00}+x_{01}&x_{00}x_{01}\\1&x_{10}+x_{20}&x_{10}x_{20}\\1&x_{11}+x_{21}&x_{11}x_{21}\end{array}\rb|\lb|\begin{array}{ccc}1&x_{00}+x_{01}&x_{00}x_{01}\\1&x_{10}+x_{20}&x_{10}x_{20}\\1&x_{11}+x_{21}&x_{11}x_{21}\end{array}\rb|=0.\label{eq}
\end{equation}
Now, after some computation, we see that the LHS of (\ref{eq}) is
\[\frac{48\sqrt{-3}\alpha_0(N^3-2N+\alpha_0(2N^2-1))}{\pm(N-1)^2},\] 
the sign depending on the choice of $\alpha_1,\alpha_2$. If $N=-1$, then, rearranging the $\alpha_i$'s if necessary, we may assume that $\alpha_0\neq-1$. We deduce that the LHS of (\ref{eq}) is non-zero, giving the required impossibility.      
\bigskip\\
We have now shown that $s_0^2-2t_0^2=1$. Define $\eta$ to be the endomorphism
\[\frac{1+s_0}{2}[1]+\frac{t_0}{2}\varepsilon+(s_0-2t_0)\psi+(-s_0+t_0)\psi\varepsilon.\]
(the coefficients are integers since $s_0^2-2t_0^2=1$ and $s_0,t_0$ are coprime). Considering the differential matrices and the images under the Rosati involution of $\varepsilon$ and $\psi$, we deduce that $\eta$ satisfies the properties described in the first paragraph of the theorem, by Proposition~\ref{dmprop} and standard properties of the Rosati involution. We have now established the first paragraph of the theorem.
\bigskip\\
Consider the second paragraph of the theorem. Suppose, for a contradiction, that the ring of endomorphisms defined over $\overline{K}$ is not generated by the endomorphisms $\varepsilon$ and $\eta$. Since the ring generated by $\varepsilon$ and $\eta$ is a maximal order in the division quaternion algebra $\lb(\frac{2,-3}{\Q}\rb)$, and the ring of endomorphisms defined over $\Qbar$ of the jacobian of $C$ is a free abelian group of finite rank (see~\cite{mil}, Proposition 12.5), the algebra of endomorphisms defined over $\Qbar$ of the jacobian of $C$ is generated by $\varepsilon,\eta,\delta$ where $\delta$ is an endomorphism lying in the centre whose square is the multiplication by $D$ endomorphism for some negative integer $D$ (see~\cite{LB}, Ch. 5 Proposition 5.7).

The strategy will be to apply Proposition~\ref{split2} and some of the theory of elliptic curves with complex multiplication. To start with, we will define some useful notation. Let $F_6(X)=\prod_{i=0}^2G_i(X)$, and let $x_{i0},x_{i1}$ be the zeros of $G_i(X)$ in $\overline{\Q}$ ($G_i(X)$ being as in (\ref{expandcurve})). Without loss of generality, we may assume that $\alpha_0$ is a real number (a rational number if $N$ is a rational cube). Let $r$ be the discriminant of $G_0$; we have
\begin{equation}
r=4(\alpha_0^2+(\alpha_0^2-1)^2).\label{disc}
\end{equation}
Since the endomorphisms $\varepsilon,\eta,\delta$ are all defined over $\Q(\sqrt{-3},\sqrt{D})$, the jacobian of $C$ splits over $\Q(\sqrt{-3},\sqrt{D})$. 
\bigskip\\
First, we will apply Proposition~\ref{split2} to show that $N$ is a rational cube, and that $\sqrt{D}$ lies in $\Q(\sqrt{-3r})$. 

The first step is to show that $r$ is positive and not a square in $\Q(\alpha_0)$. The first assertion follows immediately from (\ref{disc}). If $N$ is a rational cube, then $r$ is not a square in $\Q(\alpha_0)$; this follows from the facts that the curve $Y^2=4(X^2+(X^2-1)^2)$ is isomorphic over $\Q$ to the elliptic curve $Y^2=X^3-X^2-4X+4$ ($24A1$ in~\cite{cremona}) whose Mordell-Weil group over $\Q$ is of order $8$, and that $N\neq\pm1$ (by hypothesis). If $N$ is not a rational cube and $r$ is a square in $\Q(\alpha_0)$, equal to $(u+v\alpha_0+w\alpha_0^2)^2$ for some rational numbers $u,v,w$, then, separating the coefficients of $1,\alpha_0,\alpha_0^2$, we see, after some algebraic manipulation, that
\[v\neq0\neq w,((3w+1)v^2-(w-1)(2w+1))((3w-1)v^2+(w+1)(2w-1))=0;\]
this is impossible, again via the elliptic curve $Y^2=X^3-X^2-4X+4$.

The second step is to show that the Galois group of $F_6$ over $\Q(\alpha_0)$ is generated by $(00\:01)(10\:11)$ and $(10\:20)(11\:21)$. Since $r$ is positive and not a square in $\Q(\alpha_0)$, there is an automorphism $\sigma_1$ which alters $\sqrt{r}$ and fixes $\sqrt{-3}$, and an automorphism $\sigma_2$ which alters $\sqrt{-3}$ and fixes $\sqrt{r}$. Since the discriminant of $F_6$ is a rational square, $\sigma_1=(00\:01)(10\:11)$ (swapping $\alpha_1,\alpha_2$ if necessary), and $\sigma_2=(10\:20)(11\:21)$ (swapping $x_{20},x_{21}$ if necessary). So the Galois group of $F_6$ over $\Q(\alpha_0)$ is generated by $(00\:01)(10\:11)$ and $(10\:20)(11\:21)$, since the discriminant of $F_6$ is a rational square. 

The third step is to show that $\sqrt{D}$ lies in $\Q(\alpha_0,\sqrt{-3r})$; since $r$ is positive and not a square in $\Q(\alpha_0)$, whereas $D$ is negative, it suffices to show that $\sqrt{D}$ lies in $\Q(\alpha_0,\sqrt{-3},\sqrt{r})$ but does not lie in $\Q(\alpha_0,\sqrt{-3})$. Suppose, for a contradiction, that this is false. Then, by the second step, the Galois group of $F_6$ over $\Q(\alpha_0,\sqrt{-3},\sqrt{D})$ is generated by $(00\:01)(10\:11)$ and $(00\:01)(20\:21)$. Recall that the jacobian of $C$ splits over $\Q(\sqrt{-3},\sqrt{D})$; our contradiction follows from Proposition~\ref{split2}.

The final step is to show that $N$ is a rational cube. Suppose, for a contradiction, that this is false.  Then, $D=-3r(u+v\alpha_0+w\alpha_0^2)^2$ for some rational numbers $u,v,w$. >From (\ref{disc}), the product of the $\Q(\alpha_0,\alpha_1,\alpha_2):\Q$-conjugates of $r$ is a rational square. We deduce that $\sqrt{D}$ is $-3$ times a rational square, giving our contradiction.
\bigskip\\     
We will now apply some of the theory of elliptic curves with complex multiplication. Since the jacobian of $C$ splits over $\Q(\sqrt{-3},\sqrt{D})$ and $\sqrt{D}$ lies in $\Q(\sqrt{-3r})$, there is an abelian subvariety $E$ defined over $\Q(\sqrt{-3},\sqrt{r})$ which is an elliptic curve. This implies that
\begin{equation}
j(E)\in\Q(\sqrt{-3},\sqrt{r})\label{jinv} 
\end{equation}
and that the ring of endomorphisms of $E$ defined over $\Qbar$ is an order $\mathfrak{o}$ in $\Q(\sqrt{-3r})$; it is a standard result that
\begin{equation}
[\Q(\sqrt{-3r})(j(E)):\Q(\sqrt{-3r})]=[\Q(j(E)):\Q]=h_{\mathfrak{o}}\label{h},
\end{equation}
(see~\cite{bending}, p. 132) where $h_\mathfrak{o}$ is the class number of $\mathfrak{o}$. From (\ref{jinv}) and (\ref{h}), we see that $h_\mathfrak{o}$ divides $2$. We deduce that $r\Q^2$ lies in the set
\[\{\Q^2,2\Q^2,3\Q^2,5\Q^2,6\Q^2,17\Q^2,21\Q^2,33\Q^2,41\Q^2,57\Q^2,89\Q^2,129\Q^2,201\Q^2,489\Q^2\}\]
(see~\cite{bending}, p. 132). 
Now $r$ is not a rational square, so there is an integer $m$ not divisible by $3$ such that $r$ lies in $3m\Q^2$ or $r$ lies in $(3m+2)\Q^2$. From (\ref{disc}), rational $p,q$ solve the equation
\[q^2=t(1-p^2+p^4),\]
where $t=3m$ or $t=3m+2$ according to whether $r$ lies in $3m\Q^2$ or $r$ lies in $(3m+2)\Q^2$. However, this equation is insoluble over $\Q_3$, so it is insoluble over $\Q$, giving the required contradiction. We have now established the second paragraph of the theorem.
\finpf   
\section{Other Isogenies between Curves with $\protect\sqrt2$ Multiplication}\label{other}
\S~\ref{table} contains a table of curves defined over $\Q$ whose jacobians have $\sqrt2$ multiplication defined over $\Q$ and low conductor, which contains many pairs of curves whose jacobians are isogenous over $\Q$. Some pairs come from Theorem~\ref{isogthm} in the sense that there exist $\Delta,U,V,W\in\Q$ satisfying the conditions of the theorem such that the first curve (resp. the second curve) is isomorphic over $\Q$ to $C_1$ (resp. $C_2$). The jacobians of the curves in the pairs $\{256A_1,256A_2\}$ are isogenous over $\Q$ by Theorem~\ref{isogthm2}. The jacobians of the curves in the pairs $\{675A_1,675A_2\},\{675B_1,675B_2\}$ are isogenous over $\Q$ by Theorem~\ref{quatthm}. Taking into account that `isogenous over $\Q$' is an equivalence relation, this leaves $42$ isogenies to be discovered, corresponding to the following pairs of curves (for brevity, $29A{1,2}$ denotes $\{29A1,29A2\}$, etc.):
\begin{center}
\begin{tabular}{cccccccc}
$\!29A1,2\!$&$\!261A1,2\!$&$\!464A1,2\!$&$\!725A1,2\!$&$\!841A1,2\!$&$\!1421A1,2\!$&$\!1856A1,2\!$&$\!1856B1,2\!$\\
$\!39A1,5\!$&$\!117A1,5\!$&$\!507A1,5\!$&$\!624A1,5\!$&$\!975A1,5\!$&$\!1521A1,5\!$&$\!1872A1,5\!$&$\!1911A1,5\!$\\
$\!43A1,2\!$&$\!387A1,2\!$&$\!688A1,2\!$&$\!1075A1,2\!$&$\!1849A1,2\!$\\
$\!65A1,3\!$&$\!325A1,3\!$&$\!585A1,3\!$&$\!845A1,3\!$&$\!1040A1,3\!$\\
$\!82A^*1,3\!$&$\!656A^*1,3\!$&$\!738A^*1,3\!$\\
$\!91A1,2\!$&$\!637A1,2\!$&$\!819A1,2\!$&$\!1183A1,2\!$&$\!1456A1,2\!$\\
$\!98A2,1\!$&$\!98A4,1\!$&$\!784A2,1\!$&$\!784A4,1\!$&$\!882A2,1\!$&$\!882A4,1\!$\\
$\!529A2,1\!$&$\!529B2,1\!$
\end{tabular}
\end{center}
In fact, we can reduce the number of isogenies to be discovered from $42$ to $8$, as we now explain. To do this, we will need the following definition. If $F$ is a rational sextic with six distinct roots (possibly including $\infty$), and $\Delta$ is a non-zero rational number, then the {\em twist} by $\Delta$ of the curve $Y^2=F(X)$ is the curve $Y^2=\Delta F(X)$. If the jacobians of two curves of this form are isogenous over $\Q$, then the jacobians of their twists by $\Delta$ are, since the differential matrix of the induced isomorphism between the jacobian of a curve and the jacobian of its twist by $\Delta$ is a scalar matrix. So, if $\sim$ is the equivalence relation on the $42$ pairs of curves mentioned above defined by $\{C_1,C_2\}\sim\{C_3,C_4\}$ if and only if $C_3$ (resp. $C_4$) is isomorphic over $\Q$ to the twist by $\Delta$ of $C_1$ (resp. $C_2$) for some non-zero rational number $\Delta$, then we can reduce the number of isogenies to be discovered from $42$ to the number of equivalence classes, which is $8$; in the list above, the equivalence classes correspond to the rows, which can be seen by inspection of our table, apart from the fact that $98A1$ is isomorphic over $\Q$ to its own twist by $-7$, an isomorphism being $(x,y)\mapsto\lb(\frac{3-2x}{x+2},\frac{49y}{(x+2)^3}\rb)$. Note that any twist will do, not necessarily one that appears in our table; we derived an isogeny defined over $\Q$ between the jacobians of the twists of $98A2$ and $98A1$ by $3$, and the twists of $529A2$ and $529A1$ by $3$, none of which appear in our table (why will become apparent later). 

In this rest of this section, we will give the isogenies explicitly, and explain how they were derived. 
\bigskip\\   
Let $C_1,C_2$ be curves of genus $2$ defined over $\C$ whose jacobians appear to be homomorphic. Our method for deriving the homomorphism between the jacobians of $C_1,C_2$ uses the complex-analytic theory of abelian varieties. It was used by P. Van Wamelen to derive explicit isogenies defined over $\Q$ between jacobians of curves of genus $2$ defined over $\Q$ with good reduction away from $2$, the curves originally being discovered by N. P. Smart (see~\cite{smart}). 

First, we identify the jacobians with complex tori as described at the end of \S~\ref{isog}, choosing the symplectic bases in advance and finding the period matrices via numerical integration. Then, we derive the homomorphism $\varphi$ between the complex tori by finding the
matrices $A_\phi$ and $R_\phi$ (defined in \S~\ref{isog}); how we found
them will be explained later. Finally, we derive the homomorphism $\phi$ between the jacobians, as follows. We choose a point $P_1$ on $C_1$, and, for several points $P=(x,y)$ on $C_1$, we work out the image of 
\[[(z_1,t_1)-(z_2,-t_2)]=\phi([P-P_1])\]
on the Kummer surface for $C_2$ via $\varphi$ and the Abel-Jacobi isomorphisms, from which we obtain $z_1+z_2,z_1z_2$ in terms of $x,y$ (see~\cite{CF}, Ch. 3 for an explicit description of the Kummer surface). Eventually, we can work out what $R_1(x),\dots,R_4(x)$ are, and then we can work out what $S_1(x),\dots,S_4(x)$ are (using notation from \S~\ref{isog}), giving $\phi$ explicitly. If it is an isogeny, then its degree is the determinant of $R_\phi$.

The inverse of the Abel-Jacobi isomorphism (in particular, its image on the Kummer surface) can be expressed using theta functions, which are functions from $\C^2$ to $\C$ involving a convergent sum of terms indexed by integers $n_1,n_2$, the terms tending to zero as $n_1,n_2\rightarrow\infty$. For more details about the image on the Kummer surface, see~\cite{LB}, Ch. 10 \S 3. The Kummer surfaces in~\cite{CF} and~\cite{LB} are not necessarily equal, but nevertheless differ by an automorphism of $\p^3$, which can be determined via its action on singular points. 
\bigskip\\
For each of our isogenies, evaluating the numerical integrations to 100 significant figures and evaluating the sum of the terms with $|n_1|,|n_2|$ less than or equal to $15$ for the theta functions was more than sufficient. 
\bigskip\\          
For each isogeny in~\cite{wamelen}, $R_1(x),\dots,R_4(x),S_1(x),\dots,S_4(x)$ are defined over $\Q$, its field of definition, $R_2,R_4,S_1,S_3$ are zero, and $P_1$ is a Weierstrass point, not necessarily defined over $\Q$. For each of our isogenies, $R_1(x),\dots,R_4(x),S_1(x),\dots,S_4(x)$ are defined over $\Q$, its field of definition, $R_2,R_4,S_1,S_3$ are not necessarily zero, and $P_1$ is $\infty^+$. 

We derived an isogeny between the jacobians of the twists of the curves $98A2$ and $98A1$ by $3$, and the twists of the curves $529A2$ and $529A1$ by $3$, none of which appear in our table, to ensure that $P_1$ could be taken to be a point defined over $\Q$.
 
For the isogeny between the jacobians of the curves $29A1$ and $29A2$, $R_2,R_4,S_1,S_3$ are zero even though $P_1$ is not a Weierstrass point (the kernel of the isogeny contains $[\infty^+-\infty^-]$). A similar observation applies to three of the other pairs of curves, namely
\[\{39A1,39A5\},\{65A1,65A3\},\{82A^*1,82A^*3\}.\]
For each isogeny $\phi$ in~\cite{wamelen}, an LLL based algorithm is used to find $A_\phi$ and $R_\phi$, no restriction being placed on the endomorphism rings of the jacobians. For each of our isogenies $\phi$, finding $A_\phi$ and $R_\phi$ is considerably easier due to the existence of $\sqrt2$ multiplications, as we now explain. 

Let $C_1,C_2$ be the curves corresponding to $\phi$, let $\varepsilon_i$ be a $\sqrt2$ multiplication defined over $\Q$ of the jacobian of $C_i$ which is fixed by the Rosati involution, let $\lambda_i$ be the canonical polarisation on the jacobian of $C_i$, and let $H_i$ be the corresponding hermitian form on $\C^2$ as defined in \S~\ref{isog} (recall that we computed a corresponding period matrix when working out $\phi$ explicitly, and that the form can be computed from the matrix and is independent of the choice of the symplectic basis). We will explain how $A_\phi$ is found; then, $R_\phi$, is easily found (equation (\ref{pereqn}) ensures that the elements of $R_\phi$, which are integers, kill linear polynomials defined over $\C$; the elements are easily found from the polynomials, by trial and error using small integers). The ring of endomorphisms defined over $\Q$ of the jacobian of $C_1$ is $\Z(\sqrt2)$ (the reason for this is discussed in \S~\ref{table}), so the pullback of $\lambda_2$ by $\phi$ is $\lambda_1(n_1+n_2\varepsilon_1)$ for some integers $n_1,n_2$.  Complex-analytically, this means that
\[H_1((n_1I_2+n_2A_{\varepsilon_1})\underline{z},\underline{w})=H_2(A_\phi\underline{z},A_\phi\underline{w})\]
for all $\underline{z},\underline{w}\in\C^2$; so, letting $M_i$ be the positive definite hermitian matrix such that
\begin{equation}
H_i(\underline{z},\underline{w})=\underline{z}^tM_i\overline{\underline{w}}\label{hm}
\end{equation}
for all $\underline{z},\underline{w}\in\C^2$, we have
\begin{equation}
(n_1I_2+n_2A_{\varepsilon_1}^t)M_1=A_\phi^tM_2\overline{A_\phi}=A_\phi^tM_2A_\phi\label{comp}
\end{equation}
(the second equality follows since $\phi$ is defined over $\Q$). An intelligent guess of the degree of $\phi$, i.e. of $n_1^2-2n_2^2$ (which is $7$ in most cases) leads to an intelligent guess of $n_1,n_2$, and $A_{\varepsilon_1}$ can be easily found using Proposition~\ref{dmprop} by finding an isomorphism defined over $\Q$ from $C_1$ to its Richelot dual. Then, (\ref{comp}) is equivalent to four homogeneous quadratics defined over $\C$ in the elements of $A_\phi$, two with three terms and two with four terms, being zero; the elements are easily found from the quadratics, by trial and error using rational numbers with small numerators and denominators. 

For each of our isogenies $\phi$, $A_\phi$ is invertible. This implies that $\phi$ is indeed an isogeny, without needing to know $M_1,M_2,n_1,n_2$; indeed, $n_1,n_2$ are not both zero by (\ref{comp}), since $M_1,M_2$ are invertible, and if $D$ is the divisor on the product of $C_1$ with $C_2$ induced by $\phi$, which implies that $\phi=\phi_{[D]}$, then
\begin{equation}
\phi_{[D]}\phi_{\tau^*([D])}=n_1+n_2\varepsilon,\label{comp2}
\end{equation}
by Proposition~\ref{proprich}.    

For each of our isogenies $\phi$, $n_1$ turns out to be the smallest positive integer possible (the pullback of $\lambda_2$ by $\phi$ is $\lambda_1(n_1+n_2\varepsilon_1)$, so the zeros of the characteristic polynomial of the differential matrix of $n_1+n_2\varepsilon_1$ must be positive (see~\cite{LB}, Ch. 5 Theorem 2.4), so $n_1$ must be positive). For example, if the degree of $\phi$ is $7$, then $n_1$ is $3$. Since $n_1$ is the smallest integer possible, it is reasonable to hope that the explicit expression for $n_1+n_2\varepsilon$ is as nice as possible, so it is reasonable to hope that the explicit expression for $\phi$ is as nice as possible, from (\ref{comp2}).
\bigskip\\
The results of our computations are shown below. For each of our isogenies, we give $C_1,C_2$, the isogeny, its degree, its kernel, $M_1,M_2,A_{\varepsilon_1},A_\phi,n_1,n_2$. Before the results are shown, we make some points about the isogeny, $M_1,M_2,n_1,n_2$.

The isogeny is given explicity as follows:
\[[(x,y)-\infty^+]\mapsto[(z_1,t_1)-(z_2,-t_2)]\]
where $z_1,z_2$ are the zeros of the expression in $x,y,z$ w.r.t. $z$ and $t_k$ is given by the expression below. 

Since $\varepsilon_i$ is fixed by the Rosati involution, the pullback of the canonical polarisation on the jacobian of $C_i$ by $\varepsilon_i$ is twice the canonical polarisation, so
\[H_i(2\underline{z},\underline{w})=H_i(A_{\varepsilon_i}\underline{z},A_{\varepsilon_i}\underline{w})\]
for all $\underline{z},\underline{w}\in\C^2$, which implies, by (\ref{hm}), that
\[2M_i=A_{\varepsilon_i}^tM_i\overline{A_{\varepsilon_i}}=A_{\varepsilon_i}^tM_iA_{\varepsilon_i}\]
(the second equality follows since $\varepsilon_i$ is defined over $\Q$). We deduce that $M_i$ is a real symmetric matrix; we give its elements to eight significant figures.

For each of our isogenies, $n_2$ is positive; the sign of $\varepsilon_1$ can be changed if necessary to ensure this. 
%
\begin{eqnarray*}
C_1:Y^2\tabeq X^6-2X^5+7X^4-6X^3+13X^2-4X+8,\\
C_2:T^2\tabeq Z^6-4Z^5-12Z^4+2Z^3+8Z^2+8Z-7: 
\end{eqnarray*} 
\begin{eqnarray*}
\tabspace z^2+z(x^3+x+1)+x^4-x^3+3x^2-x+1,\\ 
t_k\tabeq y(x(x^3+x+2)+z_k(x^3+x^2-x+1)). 
\end{eqnarray*}
Degree is $7$, kernel is generated by $[\infty^+-\infty^-]$,
\begin{center}
$M_1=\lb(\begin{array}{cc}0.24487048&-0.029342177\\-0.029342177&0.12243524\end{array}\rb),M_2=\lb(\begin{array}{cc}0.52414966&-0.063750885\\-0.063750885&0.39664789\end{array}\rb),$

$A_{\varepsilon_1}=\lb(\begin{array}{cc}0&-1\\-2&0\end{array}\rb),A_\phi=\lb(\begin{array}{cc}-1&0\\-1&1\end{array}\rb),n_1=3,n_2=1.$
\bigskip\\
****************************************************************
\end{center}
\begin{eqnarray*}
C_1:Y^2\tabeq(X-1)(X+3)(X^2-X+1)(X^2-X-3),\\
C_2:T^2\tabeq-3(Z-2)(5Z+2)(Z^2-3Z-1)(Z^2+Z+7): 
\end{eqnarray*} 
\begin{eqnarray*}
\tabspace z^2(2x^2-2x+3)-3z(x^3+x^2-2x+3)+3x^4-5x^2+2x-3,\\
t_k\tabeq\frac{9xy[x(x^2-5x+1)-3z_k(x^2-2x+2)]}{(2x^2-2x+3)^2}.
\end{eqnarray*}
Degree is $7$, kernel is generated by $[\infty^+-\infty^-]$,
\begin{center}
$M_1=\lb(\begin{array}{cc}0.52150456&-0.12359786\\-0.12359786&0.27430883\end{array}\rb),M_2=\lb(\begin{array}{cc}12.099879&-0.86836281\\-0.86836281&6.0499397\end{array}\rb),$
$A_{\varepsilon_1}=\lb(\begin{array}{cc}1&-1\\-1&-1\end{array}\rb),A_\phi=\lb(\begin{array}{cc}-\frac{1}{3}&0\\\frac{1}{3}&-\frac{1}{3}\end{array}\rb),n_1=3,n_2=1.$
\bigskip\\
****************************************************************
\end{center} 
\begin{eqnarray*}
C_1:Y^2\tabeq(X-2)(X+2)(X^4+2X^3+7X^2+2X+1),\\
C_2:T^2\tabeq-3(5Z^2-12Z+8)(17Z^4+18Z^3-29Z^2-18Z+17):
\end{eqnarray*}
\begin{eqnarray*}
\tabspace z^2(3x^3+x-1)-3zx^2(2x+1)+3x^3+3x^2+x+2,\\
t_k\tabeq\frac{9y[x(6x^4+15x^3+5x^2+13x+3)-z_k(6x^5+12x^4-8x^3+5x^2+2x+1)]}{(x+2)(3x^3+x-1)^2}.
\end{eqnarray*}
Degree is $7$, kernel is generated by $[\infty^+-\infty^-]$,
\begin{center}
$M_1=\lb(\begin{array}{cc}0.71888934&0.11632882\\0.11632882&0.35944467\end{array}\rb),M_2=\lb(\begin{array}{cc}17.316093&-13.986967\\-13.986967&19.315888\end{array}\rb),$
$A_{\varepsilon_1}=\lb(\begin{array}{cc}0&-1\\-2&0\end{array}\rb),A_\phi=\lb(\begin{array}{cc}-\frac{1}{3}&\frac{1}{3}\\0&\frac{1}{3}\end{array}\rb),n_1=3,n_2=1.$
\bigskip\\
****************************************************************
\end{center} 
\begin{eqnarray*}
C_1:Y^2\tabeq(X+1)(2X+3)(2X^4-3X^3-X^2+8X-4),\\ 
C_2:T^2\tabeq-3(2Z^2-7Z+1)(2Z^4-11Z^3+55Z^2-52Z+2): 
\end{eqnarray*}
\begin{eqnarray*}
\tabspace2z^2x-z(6x^2+7x-6)+6x^3+12x^2+x-6,\\
t_k\tabeq\frac{9y[(x+1)(2x^2-3)+z_k(x+3)]}{(2x+3)x^2}.
\end{eqnarray*}
Degree is $7$, kernel is generated by $[\infty^+-\infty^-]$,
\begin{center}
$M_1=\lb(\begin{array}{cc}0.75163788&0.40201379\\0.40201379&0.42820864\end{array}\rb),M_2=\lb(\begin{array}{cc}14.944004&-4.5611388\\-4.5611388&7.4720019\end{array}\rb),$
$A_{\varepsilon_1}=\lb(\begin{array}{cc}2&1\\-2&-2\end{array}\rb),A_\phi=\lb(\begin{array}{cc}\frac{1}{3}&0\\\frac{2}{3}&\frac{1}{3}\end{array}\rb),n_1=3,n_2=1.$
\bigskip\\
****************************************************************
\end{center} 
\begin{eqnarray*}
C_1:Y^2\tabeq X^6+8X^3-36X^2+12X-4,\\
C_2:T^2\tabeq3(3Z^6+2Z^5-7Z^4+4Z^3+13Z^2-10Z+7):
\end{eqnarray*} 
\begin{eqnarray*}
\tabspace2z^2(x^2-x-1)(2x^4+5x^3+16x^2-5x+2)\\
\tabminus z(x^2-x+2)(x^5+6x^4+25x^3+12x^2-4x+2-y(x^2+2x-1))\\
\tabminus x^8-3x^7-9x^6-8x^5-18x^4+21x^3+34x^2-12x+6\\
\tabplus y(x^5+3x^4+13x^3+6x^2+8x-1),\\
t_k\tabeq2^{-1}((x^2-x-1)(2x^4+5x^3+16x^2-5x+2))^{-2}\\
\tabtimes[5x^{15}+28x^{14}+127x^{13}+278x^{12}+538x^{11}+226x^{10}+106x^9-181x^8\\
\tabminus4339x^7+736x^6+3413x^5-4945x^4+2036x^3-854x^2+142x-36\\
\tabminus y(5x^{12}+28x^{11}+127x^{10}+250x^9+412x^8+510x^7\\
\tabplus866x^6+981x^5+49x^4-82x^3+151x^2-19x+2)\\
\tabplus z_k(7x^{14}+41x^{13}+200x^{12}+364x^{11}+827x^{10}+1343x^9+143x^8\\
\tabplus2380x^7-722x^6+3212x^5-2285x^4+838x^3-206x^2+14x-4\\
\tabminus y(7x^{11}+41x^{10}+168x^9+408x^8+373x^7-815x^6\\
\tabminus905x^5-182x^4+332x^3-202x^2+57x-10))]. 
\end{eqnarray*}
Degree is $7$, kernel is generated by $[(-1+\sqrt2,3\sqrt{-3}(1-\sqrt2))-(-1-\sqrt2,-3\sqrt{-3}(1+\sqrt2))]$,
\begin{center}
$M_1=\lb(\begin{array}{cc}0.85231579&-0.10074000\\-0.10074000&0.42615790\end{array}\rb),M_2=\lb(\begin{array}{cc}1.8285695&-0.22467789\\-0.22467789&1.3792137\end{array}\rb),$
$A_{\varepsilon_1}=\lb(\begin{array}{cc}0&-1\\-2&0\end{array}\rb),A_\phi=\lb(\begin{array}{cc}1&0\\1&-1\end{array}\rb),n_1=3,n_2=1.$
\bigskip\\
****************************************************************
\end{center} 
\begin{eqnarray*}
C_1:Y^2\tabeq3(3X^6+2X^5+13X^4+12X^3+17X^2+14X+7),\\
C_2:T^2\tabeq-3(13Z^6+156Z^5+288Z^4-560Z^3-936Z^2+804Z-148):
\end{eqnarray*}
\begin{eqnarray*}
\tabspace2z^2(110x^8+115x^7+595x^6+475x^5+1173x^4+747x^3+777x^2+514x+332)\\
\tabplus z(297x^9+1243x^8+2573x^7+7196x^6+7973x^5\\
\tabplus13263x^4+10413x^3+8196x^2+4766x+1768\\
\tabplus9y(11x^6+12x^5+84x^4+23x^3+120x^2+14x+24))\\
\tabplus297x^{10}+855x^9+2861x^8+6139x^7+10471x^6\\
\tabplus16210x^5+17310x^4+17046x^3+11559x^2+4684x+698\\
\tabplus y(99x^7+252x^6+774x^5+1098x^4+1791x^3+1503x^2+1080x+234),\\
t_k\tabeq27\times2^{-1}(110x^8+115x^7+595x^6+475x^5+1173x^4+747x^3+777x^2+514x+332)^{-2}\\
\tabtimes[3(2541x^{19}+11484x^{18}+58022x^{17}+178482x^{16}+540857x^{15}\\
\tabplus1254025x^{14}+2790782x^{13}+5181731x^{12}+8999601x^{11}+13641458x^{10}\\
\tabplus18929071x^9+23163599x^8+25528969x^7+24386533x^6+20374577x^5\\
\tabplus14188396x^4+8035660x^3+3424072x^2+993868x+128400)\\
\tabplus y(2541x^{16}+10637x^{15}+49112x^{14}+141806x^{13}+365724x^{12}+809927x^{11}\\
\tabplus1504653x^{10}+2595523x^9+3773322x^8+5057040x^7+5714814x^6+5779843x^5\\
\tabplus4746844x^4+3245006x^3+1634108x^2+551116x+102416)\\
\tabminus z_k(3(1089x^{18}-3773x^{17}-16449x^{16}-73263x^{15}-270669x^{14}-599761x^{13}\\
\tabminus1485721x^{12}-2577140x^{11}-4543786x^{10}-6349927x^9-8594457x^8-9303947x^7\\
\tabminus9671042x^6-7873275x^5-5790545x^4-3294416x^3-1505942x^2-424028x-126440)\\
\tabplus y(1089x^{15}-4136x^{14}-1636x^{13}-89020x^{12}-100539x^{11}-494038x^{10}\\
\tabminus601320x^9-1371563x^8-1561229x^7-2243571x^6+1974570x^5\\
\tabminus1944979x^4+1170236x^3+641218x^2+224132x+27576))].
\end{eqnarray*}
Degree is $7$, kernel is generated by $[(\sqrt{-2},\sqrt{-3}(1+\sqrt{-2}))-(-\sqrt{-2},-\sqrt{-3}(1-\sqrt{-2}))]$,
\begin{center}
$M_1=\lb(\begin{array}{cc}1.3886222&0.19473835\\0.19473835&0.99914547\end{array}\rb),M_2=\lb(\begin{array}{cc}39.474527&3.73437378\\3.73437378&19.737264\end{array}\rb),$
$A_{\varepsilon_1}=\lb(\begin{array}{cc}1&1\\1&-1\end{array}\rb),A_\phi=\lb(\begin{array}{cc}-\frac{1}{3}&0\\\frac{1}{3}&\frac{1}{3}\end{array}\rb),n_1=3,n_2=1.$
\bigskip\\
****************************************************************
\end{center} 
\begin{eqnarray*}
C_1:Y^2\tabeq(X^2-4X-4)(4X^4-20X^3-X^2+2X-5),\\
C_2:T^2\tabeq3(Z^2+4Z-3)(20Z^4+188Z^3+221Z^2-536Z+208):
\end{eqnarray*}
\begin{eqnarray*}
\tabspace2z^2(112x^{10}-1184x^9-272x^8+30512x^7-19208x^6-286888x^5\\
\tabminus71848x^4+1091748x^3+1705679x^2+983196x+200728)\\
\tabplus z(1344x^{10}-12832x^9-7104x^8+288056x^7-46060x^6-2339988x^5\\
\tabminus1399160x^4+6593496x^3+11173221x^2+6546376x+1341056\\
\tabplus2y(112x^7-728x^6-1092x^5+10976x^4+3164x^3-51576x^2-59122x-17794))\\
\tabplus2(1120x^{10}-9104x^9-3384x^8+93516x^7+176204x^6-70728x^5\\
\tabminus1569722x^4-4258093x^3-4942212x^2-2621281x-522006\\
\tabplus y(336x^7-1456x^6-2772x^5+1176x^4+14560x^3+42084x^2+39585x+11557)),
\end{eqnarray*}
\begin{eqnarray*}
t_k\tabeq14(112x^{10}-1184x^9-272x^8+30512x^7-19208x^6-286888x^5\\
\tabminus71848x^4+1091748x^3+1705679x^2+983196x+200728)^{-2}\\
\tabtimes[25088x^{21}-482048x^{20}+2875904x^{19}-196800x^{18}-53969952x^{17}+87796128x^{16}\\
\tabplus603612864x^{15}-1228624128x^{14}-4883582592x^{13}+2764585456x^{12}\\
\tabplus37549255280x^{11}+80002576660x^{10}-113615454834x^9-778431091878x^8\\
\tabminus926460528102x^7+1053308353944x^6+4177656977325x^5+5306677525866x^4\\
\tabplus3681311854981x^3+1487688030962x^2+329899578472x+31137395904\\
\tabminus y(12544x^{18}-207872x^{17}-83968x^{16}+14358368x^{15}-31935872x^{14}-315305536x^{13}\\
\tabplus577228736x^{12}+4287236576x^{11}-1484986160x^{10}-33567857424x^9\\
\tabminus42931045840x^8+74056656622x^7+288666008596x^6+402040053310x^5\\
\tabplus314041465892x^4+146767922063x^3+39440663991x^2+5179854420x+184765904)\\
\tabplus z_k(2(6272x^{21}-139776x^{20}+891456x^{19}+1933760x^{18}-40548048x^{17}+58507344x^{16}\\
\tabplus741957792x^{15}-1951255104x^{14}-7689117576x^{13}+22421200864x^{12}\\
\tabplus57517105404x^{11}-124071984864x^{10}-331204401227x^9+262645159569x^8\\
\tabplus1215990081543x^7+488963528319x^6-1942661031369x^5-3384080139084x^4\\
\tabminus2605915054565x^3-1103888762634x^2-250572212712x-23925980864)\\
\tabminus y(6272x^{18}-106624x^{17}+215616x^{16}+4532736x^{15}-16825520x^{14}-78826944x^{13}\\
\tabplus269497648x^{12}+954809520x^{11}-1294220592x^{10}-7550256880x^9\\
\tabminus6574570396x^8+16892629812x^7+61003304621x^6+98846184822x^5\\
\tabplus99383580282x^4+64218557270x^3+25781871770x^2+5830207512x+565685600))].
\end{eqnarray*}
Degree is $7$, $[(x,y)-(x',-y')]$ lies in the kernel if and only if $(x',y')=(x,-y)$ or there exists an integer $n$ such that, defining $\alpha$ to be $2\cos(\frac{2n\pi}{7})$,
\[x+x'=4-3\alpha^2,xx'=5+12\alpha-12\alpha^2,yy'=4158+6048\alpha-7371\alpha^2,\]
\begin{center}
$M_1=\lb(\begin{array}{cc}2.3909647&0\\0&1.1954824\end{array}\rb),M_2=\lb(\begin{array}{cc}41.841883&-8.3683766\\-8.3683766&25.105130\end{array}\rb),$
$A_{\varepsilon_1}=\lb(\begin{array}{cc}0&1\\2&0\end{array}\rb),A_\phi=\lb(\begin{array}{cc}-\frac{3}{7}&\frac{1}{7}\\-\frac{1}{7}&-\frac{2}{7}\end{array}\rb),n_1=3,n_2=1.$
\bigskip\\
****************************************************************
\end{center} 
\begin{eqnarray*}
C_1:Y^2=3(3X^6-2X^5+21X^4+2X^3+111X^2-100X+24),\\  
C_2:T^2=-23(Z^6-12Z^5+24Z^4+82Z^3+84Z^2+60Z+17):
\end{eqnarray*}
\begin{eqnarray*}
\tabspace2z^2(32x^{16}-108x^{15}+44x^{14}+2836x^{13}-7045x^{12}\\
\tabminus15747x^{11}+45506x^{10}+6414x^9+125405x^8-479947x^7+338999x^6\\
\tabminus356735x^5+1029874x^4-975758x^3+397289x^2-74429x+5294)\\ 
\tabminus z(36x^{17}+68x^{16}-1224x^{15}+5948x^{14}+24907x^{13}-78286x^{12}\\
\tabminus212508x^{11}+503657x^{10}+274221x^9-308575x^8-1057666x^7-228649x^6\\
\tabplus3861148x^5-5130038x^4+2887303x^3-781336x^2+96241x-3928\\
\tabplus3y(2x^4-x^2+x-1)(2x^{10}+8x^9+19x^8-38x^7-137x^6\\
\tabminus661x^5+1267x^4-474x^3+934x^2-562x+87))\\
\tabplus36x^{18}-12x^{17}-112x^{16}+3876x^{15}+6917x^{14}-57104x^{13}+51335x^{12}\\
\tabplus129156x^{11}-173536x^{10}-156921x^9-78274x^8+1008683x^7-234670x^6\\
\tabminus1794293x^5+2229883x^4-1122422x^3+289652x^2-40760x+2756\\
\tabplus3y(2x^4-x^2+x-1)(2x^{11}-13x^9+162x^8-146x^7-228x^6\\
\tabminus91x^5-367x^4+1473x^3-1676x^2+472x-26),
\end{eqnarray*}
\begin{eqnarray*}
t_k\tabeq3\times2^{-1}(32x^{16}-108x^{15}+44x^{14}+2836x^{13}-7045x^{12}\\
\tabminus15747x^{11}+45506x^{10}+6414x^9+125405x^8-479947x^7+338999x^6\\
\tabminus356735x^5+1029874x^4-975758x^3+397289x^2-74429x+5294)^{-2}\\ 
\tabtimes[3(336x^{35}-1824x^{34}+3264x^{33}+49640x^{32}-133488x^{31}-1206596x^{30}\\
\tabplus3209890x^{29}-1943080x^{28}-21091885x^{27}+11195162x^{26}+5314283x^{25}\\ 
\tabplus492209382x^{24}-2097606629x^{23}+5384706226x^{22}-5982560935x^{21}\\
\tabplus27234320355x^{20}-118583423667x^{19}+171018630007x^{18}-186933210378x^{17}\\
\tabplus434555564876x^{16}-450337007539x^{15}-187509831395x^{14}+857794066370x^{13}\\
\tabminus1781794272139x^{12}+3535820758233x^{11}-5087328674530x^{10}+5226627198890x^9\\  
\tabminus4218269608146x^8+2783764810985x^7-1420235364429x^6+508286231885x^5\\
\tabminus112994873855x^4+11538808518x^3+633646249x^2-278047102x+19988608)\\
\tabplus y(336x^{32}-1712x^{31}+1536x^{30}+46392x^{29}-357416x^{28}+1029172x^{27}+2046812x^{26}\\
\tabminus20428980x^{25}-842575x^{24}-3660669x^{23}+140455733x^{22}+1463723036x^{21}\\
\tabminus2460377896x^{20}-11518120090x^{19}+31112806307x^{18}-21776310392x^{17}\\
\tabplus44765058148x^{16}-176798628677x^{15}+179206694085x^{14}+170836309858x^{13}\\
\tabminus387320017250x^{12}+279096610497x^{11}-851115997717x^{10}+1955385989773x^9\\
\tabminus1899527141183x^8+742228239274x^7+69928216571x^6-175645688245x^5\\
\tabplus65690743487x^4-9553692709x^3-112188563x^2+183654253x-14784874)
\end{eqnarray*}
\begin{eqnarray*}
\tabminus z_k(3(720x^{34}-3824x^{33}-5808x^{32}+184280x^{31}-58380x^{30}\\
\tabminus1970084x^{29}-3237414x^{28}+7857246x^{27}+6441690x^{26}-21770141x^{25}\\
\tabplus326908755x^{24}-140333311x^{23}+242352863x^{22}-6132586186x^{21}\\
\tabplus9676599238x^{20}-286275273x^{19}-22053225171x^{18}-7626957021x^{17}\\
\tabplus222177837563x^{16}-430388512414x^{15}+683686569943x^{14}-1702468151775x^{13}\\
\tabplus3062247652118x^{12}-3803420206815x^{11}+3969857791955x^{10}-3464326781802x^9\\
\tabplus2421157816999x^8-1479657095101x^7+859646683615x^6-434655177017x^5\\
\tabplus165665914158x^4-43413017928x^3+7294887843x^2-706112658x+29930900)\\  
\tabplus y(720x^{31}-3584x^{30}+10144x^{29}+67768x^{28}-1120348x^{27}+955624x^{26}\\
\tabplus21799772x^{25}-54600128x^{24}-132326800x^{23}+370674141x^{22}-216008226x^{21}\\
\tabplus3245089872x^{20}+46355879x^{19}-34377023991x^{18}+19416681972x^{17}\\
\tabplus81123073864x^{16}+43319971437x^{15}-111565362666x^{14}-703182585562x^{13}\\
\tabplus1301692971248x^{12}-809330078325x^{11}+1503015560674x^{10}-3274327190043x^9\\
\tabplus2811353068588x^8-423856074658x^7-1027833468244x^6+931252574123x^5\\
\tabminus403836032278x^4+104527932890x^3-16503960825x^2+1475977435x-57596040))].   
\end{eqnarray*}  
Degree is $23$, $[(x,y)-(x',-y')]$ lies in the kernel if and only if $(x',y')=(x,-y)$ or there exists an integer $n$ such that, defining $\alpha$ to be $2\cos(\frac{2n\pi}{23})$, 
\begin{eqnarray*}
x+x'\tabeq3+3\alpha-27\alpha^2-12\alpha^3+27\alpha^4+7\alpha^5-9\alpha^6-\alpha^7+\alpha^8,\\
xx'\tabeq-5-39\alpha+25\alpha^2+83\alpha^3-41\alpha^4-62\alpha^5+29\alpha^6+19\alpha^7-9\alpha^8-2\alpha^9+\alpha^{10},\\
yy'\tabeq-7038-63204\alpha-67275\alpha^2+103983\alpha^3+177813\alpha^4-47334\alpha^5\\
\tabminus140829\alpha^6+4416\alpha^7+43608\alpha^8+621\alpha^9-4554\alpha^{10}, 
\end{eqnarray*}
\begin{center}
$M_1=\lb(\begin{array}{cc}8.0401293&-1.0618817\\-1.0618817&4.0200646\end{array}\rb),M_2=\lb(\begin{array}{cc}36.787539&7.8126670\\7.8126670&21.162205\end{array}\rb),$
$A_{\varepsilon_1}=\lb(\begin{array}{cc}0&-1\\-2&0\end{array}\rb),A_\phi=\lb(\begin{array}{cc}-1&0\\1&-1\end{array}\rb),n_1=5,n_2=1.$
\end{center} 
\newpage
\section{A Table of Curves with $\protect\sqrt2$ Multiplication}\label{table}
This section contains a table of curves defined over $\Q$ whose jacobians have $\sqrt2$ multiplication defined over $\Q$ and are simple over $\Q$, and so satisfy the conditions of the generalised Shimura-Taniyama-Weil conjecture; all the curves are isomorphic over $\Q$ to a curve from Theorem~\ref{Thm3.1.1}, as expected. In each case, the existence of the $\sqrt2$ multiplication can be easily verified using Proposition~\ref{dmprop} by finding an isomorphism defined over $\Q$ from the curve to its Richelot dual, and the simplicity can be easily verified either by using Proposition~\ref{split2} where applicable, or by showing that the characteristic polynomial of the Frobenius endomorphism of the reduction at a good prime $p$ of the jacobian of the curve, which can be computed by working out the numbers of points on the curves defined over $\F_p$ and over $\F_{p^2}$ (see~\cite{CF}, Chapter 14), is irreducible (see~\cite{tate}, \S 3 Theorem 1).
\bigskip\\
For each curve, we give a polynomial $F$ of degree $5$ or $6$ such that the curve is $Y^2=F(X)$. To each curve we associate a symbol $NLi$, where $N$ is a positive integer (e.g. $29$), $L$ is either a letter (e.g. $A$) or a starred letter (e.g. $A^*$) (we will explain the difference later), and $i$ is a positive integer (e.g. $1$). In each case, the jacobian is conjecturably isogenous over $\Q$ to the abelian variety associated via the Shimura construction to a newform on $\Gamma_0(N)$. The jacobians of two of the curves are isogenous over $\Q$ if and only if they have the same $N$ and $L$. The notation `$d:j$' in the `I' column and a particular row means that the jacobian of the curve in that row (with symbol $NLi$) is $d$-isogenous over $\Q$ to the jacobian of the curve $NLj$; for example, the jacobian of the curve $29A1$ is $7$-isogenous over $\Q$ to the jacobian of the curve $29A2$.
\bigskip\\
We will now define the conductor of an abelian variety defined over $\Q$, and explain how it helps in the compilation of our table. We will first define the ramification groups of a finite Galois extension $K_2/K_1$ of $p$-adic fields, which are involved in the definition of the conductor. Let $G$ be the Galois group of $K_2/K_1$, and let $v_{K_2}$ (resp. $R_{K_2}$) be the valuation (resp. discrete valuation ring) of $K_2$. For each integer $i\geq-1$, the {\em $i^{th}$ higher ramification group} of $K_2/K_1$ is defined to be
\[G_i=\{\sigma\in G\colon v_{K_2}(\sigma(\alpha)-\alpha)\geq i+1\:\forall\alpha\in R_{K_2}\}.\]
$G_{-1}$ is the Galois group of $K_2/K_1$, the inertia group is $G_0$, and the wild inertia group is $G_1$. For each integer $i\geq0$, $G_i$ is a subgroup of $G_{i-1}$, and there is an integer $j\geq-1$ such that $G_j$ is trivial. The following lemma records a standard fact relating the ramification groups to the different, which will be used later:
\begin{lemma}\label{ram}
Let $\mkD(K_2/K_1)$ be the different of $K_2/K_1$. Then, we have
\[v_{K_2}(\mkD(K_2/K_1))=\sum_{i=0}^{\infty}(g_i-1).\]
\end{lemma}
\pf See~\cite{serre}, Ch. IV \S 1 Proposition 4.
\finpf
  
If $A$ is an abelian variety defined over a $p$-adic field $K$, then the {\em exponent of the conductor} of $A/K$ is defined to be
\begin{equation}
r_K=\dim_{\Q_l}\frac{T_l(A)\otimes_{\Z_l}\Q_l}{{(T_l(A)\otimes_{\Z_l}\Q_l)}^{I_K}}+\sum_{i=1}^\infty\frac{g_i}{g_0}\dim_{\F_l}\frac{A[l]}{{(A[l])}^{G_i}}\label{localdef}
\end{equation}
(see~\cite{LRS}), where $l$ is a prime distinct from $p$, $I_K$ is the absolute inertia group of $K$, $G_i$ is the $i^{th}$ higher ramification group of $K(A[l])/K$, and $g_i$ is the order of $G_i$; it turns out to be an integer which is independent of $l$, although this is not clear {\em a priori}. The part involving $I_K$ is called the {\em tame part}, and the part involving the $G_i$'s is called the {\em wild part}. If $A$ is an abelian variety defined over $\Q$, then the {\em conductor} of $A$ is defined to be
\[\prod_pp^{r_{\Q_p}}.\]

Let $F$ be a field of degree $g$ over $\Q$, let $N$ be a positive integer, and let $f$ be a newform on $\Gamma_1(N)$ whose field of Fourier coefficients is $F$. If $A$ is the abelian variety associated via the Shimura construction to $f$, then, by deep work due to Deligne, Langlands and Carayol, $A$ has conductor $N^g$. In particular, when $F$ is $\Q(\sqrt2)$ and $A$ is the jacobian of a curve of genus $2$ defined over $\Q$, $A$ has conductor $N^2$. So, if the jacobian of a curve of genus $2$ defined over $\Q$ with $\sqrt2$ multiplication defined over $\Q$ has conductor $N^2$, then it is conjecturably isogenous to the abelian variety associated via the Shimura construction to a newform on $\Gamma_0(N)$ whose field of Fourier coefficients is $\Q(\sqrt2)$, since two abelian varieties which are defined over $\Q$ and isogenous over $\Q$ have the same conductor (for each $p$, (\ref{localdef}) with $l$ coprime to the degree of the isogeny is the same for both); this is how the conductor of an abelian variety defined over $\Q$ helps in the compilation of our table.
\bigskip\\
We will now discuss how the conductors of the jacobians of the curves in our table were computed. Computing the conductor of the jacobian of a curve of genus $2$ defined over $\Q$ is no trivial matter in general; however, its odd part can be computed using the program {\bf genus2reduction} (\cite{liu}), so we just need to compute its exponent of $2$. We have not been completely successful with the curves in our table; the isogeny classes labelled by letters have known conductor, and the isogeny classes labelled by starred letters have conjectured conductor. Twists of curves will be important (see \S~\ref{other} for the definition of the twist). \bigskip\\
The curves whose jacobians have known conductor fall into three categories: those which have good reduction at $2$ (see~\cite{bending}, \S 2.4 for a discussion on curves of genus $2$ with good reduction at $2$), covering the isogeny classes with $N$ odd except for 
\[203A^*,259A^*,497A^*,1421A^*,1813A^*,1827A^*,\]
those whose twist by one of $-1,\pm2$ have good reduction at $2$, covering the isogeny classes
\begin{equation}
464A,624A,688A,880A,1040A,1360A,1456A,1856A,1872A,1968A,\label{goodcat}
\end{equation}
and those whose jacobians are $4$-isogenous over a quadratic field to the square of an elliptic curve, covering the isogeny classes
\begin{eqnarray}
\tabspace98A,196A,392A,392B,784A,784B,784C,784D,882A,1764A,\label{row1}\\
\tabspace160A,256A,800A,1440A,1600A,1764A\nonumber.
\end{eqnarray} 

The curves whose jacobians have conjectured conductor fall into two categories: those which have $N\leq500$, covering the isogeny classes
\begin{eqnarray}\tabspace82A^*,94A^*,172A^*,178A^*,203A^*,218A^*,226A^*,232A^*,259A^*,262A^*,\label{badcat1}\\
\tabspace276A^*,278A^*,302A^*,332A^*,334A^*,362A^*,380A^*,390A^*,424A^*,\nonumber\\
\tabspace426A^*,430A^*,434A^*,436A^*,438A^*,464A^*,488A^*,490A^*,497A^*,\nonumber
\end{eqnarray}
and those whose twist by one of $-1,\pm2$ have $2^3\nmid N\leq500$, covering the isogeny classes
\begin{eqnarray}
\tabspace656A^*,688B^*,738A^*,752A^*,828A^*,846A^*,1104A^*,\label{badcat2}\\
\tabspace1170A^*,1278A^*,1314A^*,1328A^*,1421B^*,1520A^*,1548A^*,\nonumber\\
\tabspace1602A^*,1744A^*,1813A^*,1827A^*,1900A^*,1950A^*,1962A^*.\nonumber
\end{eqnarray}
\bigskip\\
The conductor for a curve which has good reduction at $2$ is odd (a well-known consequence of good reduction at $2$).

The conductor for a curve in each isogeny class in (\ref{badcat1}) can be conjectured by finding a newform $f$ on $\Gamma_0(N)$ whose field of Fourier coefficients is $\Q(\sqrt2)$ such that, for the first few primes $p$ of good reduction for the curve, the characteristic polynomials of the Frobenius endomorphisms of the reduction at $p$ of the jacobian of the curve and of the abelian variety associated via the Shimura construction to $f$ are the same, bearing Falting's isogeny theorem (see~\cite{faltings}, \S 5 Korollar 2) in mind. The former characteristic polynomial can be computed by working out the numbers of points on the curve defined over $\F_p$ and over $\F_p^2$ (see~\cite{CF}, Chapter 14), and the latter characteristic polynomial can be computed by working out the $p^{th}$ Fourier coefficient of $f$ (this can be done using generalised versions of algorithms described in~\cite{cremona}, Ch. II for elliptic curves).     
\bigskip\\
The curves in each isogeny class in (\ref{goodcat}) and in (\ref{badcat2}) can be dealt with using the following fact: if $C$ is a curve of genus $2$ defined over $\Q$ whose jacobian has $\sqrt2$ multiplication defined over $\Q$ and conductor $N^2$ for some positive integer for which $N$ not divisible by $2^3$, then the exponent of $2$ of the conductor of the jacobian of the twist of $C$ by $-1$ (resp. $\pm2$) is $2^8$ (resp. $2^{12}$). The twists by $-1,\pm2$ can be treated in analogous ways, so we will just consider the twist by $-1$. Let $C_{-1}$ be the twist of $C$ by $-1$, and let $J$ (resp. $J_{-1}$) be the jacobian of $C$ (resp. $C_{-1}$).

The first step is to show that $\Q_2(J[3])/\Q_2$ is at worst tamely ramified (the reason for choosing the prime $3$ will become apparent). It suffices to show that the inertia group of $\Q_2(J[3])/\Q_2$ has either odd order or an automorphism negating the $3$-torsion points on $J$; then, by (\ref{localdef}), using the fact that $N$ is not divisible by $2^3$, we deduce that $\Q_2(J[3])/\Q_2$ is at worst tamely ramified. 

Let $\varepsilon$ be a $\sqrt2$ multiplication defined over $\Q$ and fixed by the Rosati involution. Since $3$ is inert in $\Q(\sqrt2)$, there is a (non-canonical) embedding $\iota$ from the Galois group of $\Q_2(J[3])/\Q_2$ to $\GL_2(\F_{3^2})$, induced by the action of that group on the $3$-torsion points on $J$ (if $n_1,n_2$ are integers, the endomorphism $[n_1]+[n_2]\varepsilon$ can be canonically regarded as an element of $\F_{3^2}$). It suffices to show that the image of the inertia group of $\Q_2(J[3])/\Q_2$ under $\iota$ is a subgroup of $\SL_2(\F_{3^2})$, since that group contains a unique element of order $2$, namely the negation matrix.  

Let $\omega$ be a primitive cube root of unity; since $\Q_2(\omega)/\Q_2$ is unramified, we may replace $\Q_2(J[3])/\Q_2$ by $\Q_2(\omega)(J[3])/\Q_2(\omega)$. Let $e_3:J[3]\times J[3]\mapsto\mu_3$ be the Weil pairing corresponding to $3$ times the canonical polarisation on $J$. Since $\varepsilon$ is fixed by the Rosati involution, $e_3(P,\varepsilon(P))=e_3(\varepsilon(P),P)$ for all points $P$, and $e_3(\varepsilon(P_1),\varepsilon(P_2))=e_3(2P_1,P_2)$ for all points $P_1,P_2$ (the pullback of the canonical polarisation by $\varepsilon$ is twice the canonical polarisation). So, since $e_3$ is bilinear, alternating and non-degenerate, there is a basis $P_1,P_2$ for the $3$-torsion points on $J$ as a $\F_{3^2}$-vector space satisfying\\
\begin{eqnarray}
\tabspace e_3(P_1,\varepsilon(P_1))=1,e_3(P_1,P_2)=\omega,e_3(P_1,\varepsilon(P_2))=1,\label{wp}\\
\tabspace e_3(P_2,\varepsilon(P_1))=1,e_3(P_2,\varepsilon(P_2))=1,e_3(\varepsilon(P_1),\varepsilon(P_2))=\omega^2.\nonumber
\end{eqnarray}
Let $\sigma$ be an automorphism in the inertia group of $\Q_2(\omega)(J[3])/\Q_2(\omega)$; we need to show that $\iota(\sigma)$ is an element of $\SL_2(\F_{3^2})$. Let $a_1,a_2,b_1,b_2,c_1,c_2,d_1,d_2$ be integers such that
\[\sigma(P_1)=([a_1]+[a_2]\varepsilon)P_1+([b_1]+[b_2]\varepsilon)P_2,\sigma(P_2)=([c_1]+[c_2]\varepsilon)P_1+([d_1]+[d_2]\varepsilon)P_2.\]
Then, from (\ref{wp}), we have
\begin{eqnarray*}
\omega\tabeq\sigma(\omega)=\sigma(e_3(P_1,P_2))=e_3(\sigma(P_1),\sigma(P_2))\\
\tabeq e_3([a_1]+[a_2]\varepsilon)P_1+([b_1]+[b_2]\varepsilon)P_2,[c_1]+[c_2]\varepsilon)P_1+([d_1]+[d_2]\varepsilon)P_2)\\
\tabeq\omega^{a_1d_1+2a_2d_2-b_1c_1-2b_2c_2},\\
1\tabeq\sigma(1)=\sigma(e_3(P_1,\varepsilon(P_2))=e_3(\sigma(P_1),\sigma\varepsilon(P_2))=e_3(\sigma(P_1),\varepsilon\sigma(P_2))\\
\tabeq e_3([a_1]+[a_2]\varepsilon)P_1+([b_1]+[b_2]\varepsilon)P_2,[2c_2]+[c_1]\varepsilon)P_1+([2d_2]+[d_1]\varepsilon)P_2)\\
\tabeq\omega^{2a_1d_2+2a_2d_1-2b_1c_2-2b_2c_1}
\end{eqnarray*}
(since $e_3$ is bilinear, alternating and Galois invariant). So $a_1d_1+2a_2d_2-b_1c_1-2b_2c_2$ is congruent to $1$ modulo $3$, $2a_1d_2+2a_2d_1-2b_1c_2-2b_2c_1$ is congruent to $0$ modulo $3$. We deduce that $\iota(\sigma)$ is an element of $\SL_2(\F_{3^2})$, as required. 

Having shown that $\Q_2(J[3])/\Q_2$ is at worst tamely ramified (let $e$ be its ramification index, which is odd), the second step is to show that the inertia group of $\Q_2(J_{-1}[3])/\Q_2$ has order $2e$, that the next $\frac{e}{2}$ higher ramification groups have order $2$, and that the rest are trivial; then, since the inertia group of $\Q_2(J_{-1}[3])/\Q_2$ has an automorphism negating the $3$-torsion points on $J_{-1}$ (by the reasoning above), the exponent of $2$ of the conductor of $J_{-1}$ is $2^8$ by (\ref{localdef}), as required.

By the tower law for differents (applied to $\Q_2(J[3])(\sqrt{-1})/\Q_2(\sqrt{-1})/\Q_2$), we have
\begin{eqnarray}
\tabspace e(\Q_2(J[3])(\sqrt{-1})/\Q_2(\sqrt{-1}))v_{\Q_2(\sqrt{-1})}\mkD(\Q_2(\sqrt{-1})/\Q_2)\label{different}\\
\tabeq v_{\Q_2(J[3])(\sqrt{-1})}\mkD(\Q_2(J[3])(\sqrt{-1})/\Q_2)-v_{\Q_2(J[3])(\sqrt{-1})}\mkD(\Q_2(J[3])(\sqrt{-1})/\Q_2(\sqrt{-1})).\nonumber
\end{eqnarray}
Now $\Q_2(J_{-1}[3])/\Q_2$ has ramification index $2e$, since $\Q_2(J[3])/\Q_2$ has ramification index $e$; indeed, $\Q_2(J[3])(\sqrt{-1})=\Q_2(J_{-1}[3])(\sqrt{-1})$, which has ramification index $2e$ over $\Q_2$, is the composite of $\Q_2(J[3])$ and $\Q_2(J_{-1}[3])$, and if two fields are tamely ramified over the same field, then their composite is as well. Also, $v_{\Q_2(\sqrt{-1})}\mkD(\Q_2(\sqrt{-1})/\Q_2)=2$ and $v_{\Q_2(J[3])(\sqrt{-1})}\mkD(\Q_2(J[3])(\sqrt{-1})/\Q_2(\sqrt{-1}))=e-1$, by Lemma~\ref{ram}. So, from (\ref{different}), we have
\[v_{\Q_2(J_{-1}[3])}\mkD(\Q_2(J_{-1}[3])/\Q_2)=3e-1.\]
We deduce, by Lemma~\ref{ram}, that the inertia group of $\Q_2(J_{-1}[3])/\Q_2$ has order $2e$, that the next $\frac{e}{2}$ higher ramification groups have order $2$, and that the rest are trivial, as required. 
\bigskip\\ 
The first curve in each isogeny class in the first row (resp. the second row) of (\ref{row1}) is isomorphic over $\Q$ to a curve from Theorem~\ref{isogthm} with $V=2$, the field of definition of $\alpha$ satisfying (\ref{alpharel}) is $\Q(\sqrt{-7})$ (resp. $\Q(\sqrt{-1})$), and the elliptic curve $E$ is (\ref{ellcur}). 

In the first case, since $\sqrt{-7}$ lies in $\Q_2$, the exponent of $2$ of the conductor is the square of the exponent of the conductor of $E/\Q_2$, which can be computed using Tate's algorithm (see~\cite{Tat2}, \S 6).      

The second case is more complicated than the first case, since $\Q_2(\sqrt{-1})/\Q_2$ is ramified, whereas $\Q_2(\sqrt{-7})/\Q_2$ is unramified. We will show that the exponent of $2$ of the conductor for the curve $160A1$ is $10$ (as expected since $2^{10}$ exactly divides $160^2$). The first curve in each of the other isogeny classes can then be treated in an analogous manner. 

Let $J$ be the jacobian of the curve $160A1$, and let $E$ be the elliptic curve
\[E\colon T^2=(14+2\sqrt{-1})(Z^3+1)-(22+58\sqrt{-1})(Z^2+Z).\]
There is a $4$-isogeny defined over $\Q(\sqrt{-1})$ from $J$ to $E$, given explicitly by
\[[(x,y)-(x',y')]\mapsto\lb((z^2,t)-({z'}^2,t'),\lb(\frac{1}{z^2},\frac{t}{z^3}\rb)-\lb(\frac{1}{{z'}^2},\frac{t'}{{z'}^3}\rb)\rb),\]
where
\[z=\frac{x-1-\sqrt{-1}}{\sqrt{-1}x-1+\sqrt{-1}},t=\frac{16(1-\sqrt{-1})y}{(\sqrt{-1}x-1+\sqrt{-1})^3},\]
$z',t'$ having analogous definitions.
 
Using notation from (\ref{localdef}), let $G_{i,E},g_{i,E}$ (resp. $G_{i,J},g_{i,J}$) be $G_i,g_i$ when $A$ is $E$ and $K$ is $\Q_2(\sqrt{-1})$ (resp. when $A$ is $J$ and $K$ is $\Q_2$), for all $i\geq0$, and let $l$ be $3$ (the reason for choosing the prime $3$ will become apparent). 

The first step is to show that
\begin{equation}
g_{0,E}=g_{1,E}.\label{allwild}
\end{equation}
From~\cite{DSS}, \S 2, we see that the field generated by the $Z$-coordinates of the $3$-torsion points on $E$ is the splitting field of the quartic $\alpha^4+3\sqrt{-1}\alpha^2+(-1+\sqrt{-1})\alpha-3$. Its resolvent cubic has a zero in $\Q_2(\sqrt{-1})$ by Hensel's Lemma, which implies that its Galois group is a $2$-group. We deduce that (\ref{allwild}) holds, as required.

The second step is to show that
\begin{equation}
2g_{0,E}=g_{0,J},2g_{1,E}=g_{1,J}.\label{fact1}
\end{equation}
Since $160A1$ is isomorphic over $\Q$ to its own twist by $-1$, an isomorphism being $(x,y)\mapsto(-x,y)$, $\Q_2(J[3])$ contains $\sqrt{-1}$, which implies that $\Q_2(\sqrt{-1})(E[3])=\Q_2(J[3])$. We deduce that (\ref{fact1}) holds, as required.

The third step is to show that
\begin{equation}
g_{i,E}=g_{i,J}\:\forall i\geq2.\label{fact2}
\end{equation}
By the tower law for differents (applied to $\Q_2(\sqrt{-1})(E[3])=\Q_2(J[3])/\Q_2(\sqrt{-1})/\Q_2$), we have
\begin{eqnarray*}
\tabspace e(\Q_2(\sqrt{-1})(E[3])/\Q_2(\sqrt{-1}))v_{\Q_2(\sqrt{-1})}\mkD(\Q_2(\sqrt{-1})/\Q_2)\label{diff3}\\
\tabeq v_{\Q_2(J[3])}\mkD(\Q_2(J[3])/\Q_2)-v_{\Q_2(\sqrt{-1})(E[3])}\mkD(\Q_2(\sqrt{-1})(E[3])/\Q_2(\sqrt{-1})),
\end{eqnarray*}
and so, by Lemma~\ref{ram}, we have
\[2g_{0,E}=\sum_{i=0}^\infty(g_{i,J}-1)-\sum_{i=0}^\infty(g_{i,E}-1).\]
Since (\ref{allwild}) and (\ref{fact1}) hold, and $g_{i,E}\leq g_{i,J}$ for all $i\geq2$, we deduce that (\ref{fact2}) holds, as required. 
 
The fourth step is to show that
\begin{equation}
\dim_{\F_3}\frac{E[3]}{(E[3])^{G_{i,E}}}=2,\dim_{\F_3}\frac{J[3]}{(J[3])^{G_{i,J}}}=4\:\forall\:0\leq i<j,\label{fact3}
\end{equation}
where $j$ is the smallest integer such that $G_{j,E}$ is trivial. There is a (non-canonical) embedding $\iota$ from the Galois group of $\Q_2(\sqrt{-1})(E[3])/\Q_2(\sqrt{-1})$ to $\GL_2(\F_3)$, induced by the action of that group on the $3$-torsion points on $E$ (if $n$ is an integer, the endomorphism $[n]$ can be canonically regarded as an element of $\F_3$). It suffices to show that the image of the inertia group of $\Q_2(\sqrt{-1})(E[3])/\Q_2(\sqrt{-1})$ under $\iota$  embedding $\iota$ from the inertia group of $\Q_2(\sqrt{-1})(E[3])/\Q_2(\sqrt{-1})$ is a subgroup of $\SL_2(\Z/3\Z)$; indeed, $2$ divides $g_{i,E}$ for all $0\leq i<j$, since (\ref{allwild}) holds and $G_{1,E}$ is a $2$-group, and $\SL_2(\Z/3\Z)$ contains a unique element of order $2$, namely the negation matrix. 

Let $\omega$ be a primitive cube root of unity; since $\Q_2/\Q_2(\omega)$ is unramified, we may replace $\Q_2(\sqrt{-1})(E[3])/\Q_2(\sqrt{-1})$ by $\Q_2(\omega)(\sqrt{-1})(E[3])/\Q_2(\omega)(\sqrt{-1})$. Let $e_3:E[3]\times E[3]\rightarrow\mu_3$ be the Weil pairing corresponding to $3$ times the canonical polarisation on $E$. Since $e_3$ is bilinear, alternating and non-degenerate, there is a basis $P_1,P_2$ for the $3$-torsion points on $E$ as a $\F_3$-vector space satisfying
\begin{equation}
e_3(P_1,P_1)=1,e_3(P_1,P_2)=\omega,e_3(P_2,P_2)=1.\label{wp3}
\end{equation}
Let $\sigma$ be an automorphism in the inertia group of $\Q_2(\omega)(\sqrt{-1})(E[3])/\Q_2(\omega)(\sqrt{-1})$; we need to show that $\iota(\sigma)$ is an element of $\SL_2(\F_3)$. Let $a,b,c,d$ be integers such that
\[\sigma(P_1)=[a]P_1+[b]P_2,\sigma(P_2)=[c]P_1+[d]P_2.\]
Then, from (\ref{wp3}), we have
\[\omega=\sigma(\omega)=\sigma(e_3(P_1,P_2))=e_3(\sigma(P_1,P_2))=e_3([a]P_1+[b]P_2,[c]P_1+[d]P_2)=\omega^{ad-bc}\]
(since $e_3$ is bilinear, alternating and Galois invariant). So $ad-bc$ is congruent to $1$ modulo $3$. We deduce that $\iota(\sigma)$ is an element of $\SL_2(\F_3)$, as required.

Since the exponent of the conductor of $E/\Q_2(\sqrt{-1})$, which can be computed using Tate's algorithm (see~\cite{tate}, \S 6), is $6$, $j\geq2$, by (\ref{localdef}), and the exponent of $2$ of the conductor for the curve $160A1$ is $10$, by (\ref{localdef}), (\ref{fact1}), (\ref{fact2}) and (\ref{fact3}), as required. 
\bigskip\\
Finally in this section, we will discuss how the curves in our table were found. Most of them were found by substituting into Theorem~\ref{Thm3.1.1}, applying~\ref{isogthm} to curves already found, and taking twists of curves already found. The search concentrated on curves with good reduction at $2$ (bearing in mind that the conductor of the jacobian of such a curve can be computed using {\bf genus2reduction}). However, it is reasonable to expect that some curves defined over $\Q$ whose jacobians have $\sqrt2$ multiplication and low conductor are hard to find by substitution. 

Our table is complete in the following sense: for each newform $f$ on $\Gamma_0(N)$ whose field of Fourier coefficients is $\Q(\sqrt2)$ ($N$ being a positive integer less than or equal to $500$), our table contains a curve of genus $2$ whose jacobian is conjecturably isogenous over $\Q$ to the abelian variety associated via the Shimura construction to $f$. For each $N$, the number of newforms can be found using~\cite{muller}. We have not shown that any of the jacobians of our curves satisfy the generalised Shimura-Taniyama-Weil conjecture; however, the jacobians of the curves in the isogeny classes $29A,39A,85A,165A$ do, since
\begin{eqnarray*}
\tabspace29A_2=X_0(29),39A_3=X_0(39)/W_{13},\\
\tabspace85A_1=X_0(85)/\langle W_5,W_{17}\rangle,165A_1=X_0(165)/\langle W_3,W_5,W_{11}\rangle
\end{eqnarray*}
(see~\cite{galbraith}, \S 4.3).

Let $F$ be a quadratic extension of $\Q$. A natural question to ask is whether a polarised abelian surface defined over $\Q$ whose ring of endomorphisms defined over $\Q$ is isomorphic to an order in $F$ is isogenous over $\Q$ to a principally polarised abelian surface defined over $\Q$ whose ring of endomorphisms defined over $\Q$ is isomorphic to a maximal order in $F$. Our table shows that the answer is yes for the abelian surfaces discussed above (the abelian variety associated via the Shimura construction to $f$ has a polarisation defined over $\Q$ induced by the canonical polarisation on the jacobian of $X_0(N)$). More generally, this question has been discussed in~\cite{wilson}, \S 4.5; if $d$ is the degree of the polarisation and $c$ is the conductor of the order, the answer is yes if $d$ is coprime to $2c$ (see~\cite{wilson}, Theorems 4.5.4 and 4.5.7); note that the converse is false.

To make sure that our table was complete, some of the curves were found using the following method. Let $f$ be a newform that was `missing' from our table. First, we worked out the $p^{th}$ Fourier coefficient of $f$ (this can be done using generalised versions of algorithms described in~\cite{cremona}, Ch. II for elliptic curves) for primes $p$ less than $10,000$. Then, we expressed the abelian variety associated via the Shimura construction complex-analytically (see~\cite{cremona}, Proposition 2.10.1), and found a principal polarisation on it, using the theory described in~\cite{LB}, Ch. 4 \S 2 (it is not clear that a principal polarisation always exists). Then, we found the absolute invariants $i_1,i_2,i_3$ of a curve of genus $2$ whose jacobian is isomorphic to it (see~\cite{wang}); they turned out to be conjecturably rational numbers (it is not clear that the absolute invariants are always conjecturably rational numbers). Then, we found a curve of genus $2$ defined over $\Q$ with those absolute invariants, by intersecting a conic with a cubic as described in~\cite{mestre}, a point on the conic being found using an algorithm in~\cite{cremona2} (it is not clear that the curve is always defined over $\Q$). In each case, the curve came from Theorem~\ref{Thm3.1.1}, as expected.       
\newpage
\newpage
\begin{center}
\begin{tabular}{|ccclc|}\hl
&$\!\!\!C\!\!\!$&&$F(X)$&$I$\\\hl
$29\!\!\!$&$\!\!\!A\!\!\!$&$\!\!\!1$&$X^6-2X^5+7X^4-6X^3+13X^2-4X+8$&$7\!:\!2$\\ 
&&$\!\!\!2$&$X^6-4X^5-12X^4+2X^3+8X^2+8X-7$&\\ 
&&&&\\
$39\!\!\!$&$\!\!\!A\!\!\!$&$\!\!\!1$&$(X-1)(X+3)(X^2-X+1)(X^2-X-3)$&$7\!:\!5$\\ 
&&$\!\!\!2$&$(X^2-X+1)(X^4+X^3+4X^2-3X+9)$&$2\!:\!1$\\ 

&&$\!\!\!3$&$(X-1)(X+3)(X^2+3X-1)(X^2-5X+3)$&$2\!:\!1$\\ 
&&$\!\!\!4$&$-5(X^2-4X-8)(X^2+X+7)(3X^2+3X+1)$&$2\!:\!1$\\
\end{tabular}
\newpage
\begin{tabular}{|ccclc|}\hl
&$\!\!\!C\!\!\!$&&$F(X)$&$I$\\\hl
$117\!\!\!$&$\!\!\!A\!\!\!$&$\!\!\!1$&$-3(X-1)(X+3)(X^2-X+1)(X^2-X-3)$&$7\!:\!5$\\ 

&&$\!\!\!2$&$-3(X^2-X+1)(X^4+X^3+4X^2-3X+9)$&$2\!:\!1$\\ 
			  
&&$\!\!\!3$&$-3(X-1)(X+3)(X^2+3X-1)(X^2-5X+3)$&$2\!:\!1$\\ 
&&$\!\!\!4$&$15(X^2-4X-8)(X^2+X+7)(3X^2+3X+1)$&$2\!:\!1$\\ 
&&$\!\!\!5$&$(X-2)(5X+2)(X^2+X+7)(X^2-3X-1)$&\\
&&$\!\!\!6$&$7(3X^2+11X+9)(9X^4-39X^3+7X^2+156X+144)$&$2\!:\!5$\\

&&$\!\!\!7$&$-5(X^2-14X-3)(X^2+X-3)(23X^2-37X-9)$&$2\!:\!5$\\
\end{tabular}
\newpage
\begin{tabular}{|ccclc|}\hl
&$\!\!\!C\!\!\!$&&$F(X)$&$I$\\\hl
$207\!\!\!$&$\!\!\!A\!\!\!$&$\!\!\!1$&$(X^2-X+1)(X^4+X^3-3X+9)$&\\        

&&$\!\!\!2$&$(X-1)(X+3)(X^4-2X^3-5X^2+6X-3)$&$2\!:\!1$\\[6pt] 
&$\!\!\!B\!\!\!$&$\!\!\!1$&$-3(X^2-X+1)(X^4+X^3-3X+9)$&\\ 
&&$\!\!\!2$&$-3(X-1)(X+3)(X^4-2X^3-5X^2+6X-3)$&$2\!:\!1$\\ 
&&&&\\                                 
$209\!\!\!$&$\!\!\!A\!\!\!$&$\!\!\!1$&$X^6+4X^5+8X^4+8X^3+8X^2-4X+4$&\\ 
&&&&\\
\end{tabular}
\newpage
\begin{tabular}{|ccclc|}\hl
&$\!\!\!C\!\!\!$&&$F(X)$&$I$\\\hl
$275\!\!\!$&$\!\!\!A\!\!\!$&$\!\!\!1$&$-5(X+1)(X-3)(X^2+X-1)(3X^2-X+1)$&\\ 

&&$\!\!\!2$&$-15(X^2-X-1)(X^4-5X^3+23X^2+20X+16)$&$2\!:\!1$\\ 

&&$\!\!\!3$&$-5(X^2+4)(X^2+X+3)(3X^2-X+1)$&$2\!:\!1$\\ 
&&$\!\!\!4$&$-15(X^2-8X-4)(X^2+X-1)(X^2-11X-1)$&$2\!:\!1$\\ 
&&&&\\
$276\!\!\!$&$\!\!\!A^*\!\!\!$&$\!\!\!1$&$-(X-1)(2X-1)(4X^4-13X^2+12)$&\\

&&$\!\!\!2$&$-2(2X^2+2X-1)(48X^4+96X^3+44X^2+4X-1)$&$2\!:\!1$\\
\end{tabular}
\newpage
\begin{tabular}{|ccclc|}\hl
&$\!\!\!C\!\!\!$&&$F(X)$&$I$\\\hl
$335\!\!\!$&$\!\!\!A\!\!\!$&$\!\!\!1$&$X^6+4X^5-48X^2+20X-4$&\\ 
&&&&\\
$339\!\!\!$&$\!\!\!A\!\!\!$&$\!\!\!1$&$-3X^6+6X^5-25X^4+36X^3-61X^2+42X-31$&\\[6pt] 

&$\!\!\!B\!\!\!$&$\!\!\!1$&$-89(9X^3-114X^2+364X+173)(87X^3+130X^2+64X+43)$&\\
&&&&\\
$345\!\!\!$&$\!\!\!A\!\!\!$&$\!\!\!1$&$X^6+12X^5+32X^4-24X^3+8X^2+12X+4$&\\ 
&&&&\\
\end{tabular}
\newpage
\begin{tabular}{|ccclc|}\hl
&$\!\!\!C\!\!\!$&&$F(X)$&$I$\\\hl
$425\!\!\!$&$\!\!\!A\!\!\!$&$\!\!\!1$&$5(X^2+4)(X^4+2X^3+3X^2-2X+1)$&\\
&&$\!\!\!2$&$5(X^2+X-1)(X^4+3X^3+6X^2-5X-1)$&$2\!:\!1$\\
&&&&\\
$426\!\!\!$&$\!\!\!A^*\!\!\!$&$\!\!\!1$&$-12X^6+12X^5+17X^4-62X^3+41X^2+12X-12$&\\
&&&&\\
$429\!\!\!$&$\!\!\!A\!\!\!$&$\!\!\!1$&$(X^2+X+3)(9X^4+21X^3+31X^2+4X+12)$&\\
&&$\!\!\!2$&$(X^2-6X-3)(X^4-18X^3-5X^2+18X+13)$&$2\!:\!1$\\
\end{tabular}
\newpage
\begin{tabular}{|ccclc|}\hl
&$\!\!\!C\!\!\!$&&$F(X)$&$I$\\\hl
$495\!\!\!$&$\!\!\!A\!\!\!$&$\!\!\!1$&$3(X+1)(X-3)(X^2+X-1)(3X^2-X+1)$&\\

&&$\!\!\!2$&$(X^2-X-1)(X^4-5X^3+23X^2+20X+16)$&$2\!:\!1$\\ 

&&$\!\!\!3$&$3(X^2+4)(X^2+X+3)(3X^2-X+1)$&$2\!:\!1$\\ 

&&$\!\!\!4$&$(X^2-8X-4)(X^2+X-1)(X^2-11X-1)$&$2\!:\!1$\\[6pt] 

&$\!\!\!B\!\!\!$&$\!\!\!1$&$-3X(X+4)(X^2+X-1)(X^2+X+3)$&\\ 

&&$\!\!\!2$&$-3(X^2-X+1)(X^4+X^3+8X^2-3X+9)$&$2\!:\!1$\\ 

&&$\!\!\!3$&$(X-2)(X+2)(X^2-3X+5)(5X^2-3X+1)$&$2\!:\!1$\\
\end{tabular}
\newpage
\begin{tabular}{|ccclc|}\hl
&$\!\!\!C\!\!\!$&&$F(X)$&$I$\\\hl
$605\!\!\!$&$\!\!\!A\!\!\!$&$\!\!\!1$&$11(X+1)(X-3)(X^2+X-1)(3X^2-X+1)$&\\
&&$\!\!\!2$&$-33(X^2-X-1)(X^4-5X^3+23X^2+20X+16)$&$2\!:\!1$\\ 
&&$\!\!\!3$&$11(X^2+4)(X^2+X+3)(3X^2-X+1)$&$2\!:\!1$\\
&&$\!\!\!4$&$-33(X^2-8X-4)(X^2+X-1)(X^2-11X-1)$&$2\!:\!1$\\ 
&&&&\\
$609\!\!\!$&$\!\!\!A\!\!\!$&$\!\!\!1$&$(X^2-X+1)(X^4+X^3+16X^2-3X+9)$&\\
&&$\!\!\!2$&$(X-1)(X+3)(X^4-2X^3+11X^2-10X-3)$&$2\!:\!1$\\
\end{tabular}
\newpage
\begin{tabular}{|ccclc|}\hl
&$\!\!\!C\!\!\!$&&$F(X)$&$I$\\\hl
$688\!\!\!$&$\!\!\!A\!\!\!$&$\!\!\!1$&$3(X^6+8X^3-36X^2+12X-4)$&$7\!:\!2$\\ 
&&$\!\!\!2$&$3X^6+2X^5-7X^4+4X^3+13X^2-10X+7$&\\ 
&$\!\!\!B^*\!\!\!$&$\!\!\!1$&$10(4X^6-4X^5-76X^4+44X^3+912X^2-2176X+1438)$&\\
&&&&\\
$703\!\!\!$&$\!\!\!A\!\!\!$&$\!\!\!1$&$X^6+2X^5-13X^4-4X^3+15X^2+14X+5$&\\ 
&&&&\\
$705\!\!\!$&$\!\!\!A\!\!\!$&$\!\!\!1$&$9X^6+12X^5-32X^4+6X^3+76X^2-60X+25$&\\
\end{tabular}
\newpage
\begin{tabular}{|ccclc|}\hl
&$\!\!\!C\!\!\!$&&$F(X)$&$I$\\\hl
$753\!\!\!$&$\!\!\!A\!\!\!$&$\!\!\!1$&$9X^6+12X^5+64X^4+40X^3+112X^2-12X+36$&\\ 
&&&&\\
$765\!\!\!$&$\!\!\!A\!\!\!$&$\!\!\!1$&$-3(X^2+4)(X^4+2X^3+3X^2-2X+1)$&\\
&&$\!\!\!2$&$-3(X^2+X-1)(X^4+3X^3+6X^2-5X-1)$&$2\!:\!1$\\
&&&&\\
$767\!\!\!$&$\!\!\!A\!\!\!$&$\!\!\!1$&$X^6+4X^5+16X^4+34X^3+60X^2+52X+17$&\\ 
&&&&\\
\end{tabular}
\newpage
\begin{tabular}{|ccclc|}\hl
&$\!\!\!C\!\!\!$&&$F(X)$&$I$\\\hl
$828\!\!\!$&$\!\!\!A^*\!\!\!$&$\!\!\!1$&$3(X-1)(2X-1)(4X^4-13X^2+12)$&\\
&&$\!\!\!2$&$6(2X^2+2X-1)(48X^4+96X^3+44X^2+4X-1)$&$2\!:\!1$\\
&&&&\\
$837\!\!\!$&$\!\!\!A\!\!\!$&$\!\!\!1$&$(X^3+6X^2+6X-3)(5X^3+6X^2-6X+1)$&\\[6pt] 
&$\!\!\!B\!\!\!$&$\!\!\!1$&$-3(X^3+6X^2+6X-3)(5X^3+6X^2-6X+1)$&\\ 
&&&&\\
$841\!\!\!$&$\!\!\!A\!\!\!$&$\!\!\!1$&$29(X^6-2X^5+7X^4-6X^3+13X^2-4X+8)$&$7\!:\!2$\\
\end{tabular}
\newpage
\begin{tabular}{|ccclc|}\hl
&$\!\!\!C\!\!\!$&&$F(X)$&$I$\\\hl
$975\!\!\!$&$\!\!\!A\!\!\!$&$\!\!\!1$&$5(X-1)(X+3)(X^2-X+1)(X^2-X-3)$&$7\!:\!5$\\ 
&&$\!\!\!2$&$5(X^2-X+1)(X^4+X^3+4X^2-3X+9)$&$2\!:\!1$\\ 

&&$\!\!\!3$&$5(X-1)(X+3)(X^2+3X-1)(X^2-5X+3)$&$2\!:\!1$\\ 
&&$\!\!\!4$&$-(X^2-4X-8)(X^2+X+7)(3X^2+3X+1)$&$2\!:\!1$\\ 
&&$\!\!\!5$&$-15(X-2)(5X+2)(X^2+X+7)(X^2-3X-1)$&\\
&&$\!\!\!6$&$-105(3X^2+11X+9)(9X^4-39X^3+7X^2+156X+144)$&$2\!:\!5$\\
&&$\!\!\!7$&$-3(X^2-14X-3)(X^2+X-3)(23X^2-37X-9)$&$2\!:\!5$\\
\end{tabular}
\newpage
\begin{tabular}{|ccclc|}\hl
&$\!\!\!C\!\!\!$&&$F(X)$&$I$\\\hl
$1183\!\!\!$&$\!\!\!A\!\!\!$&$\!\!\!1$&$-13(3X^6+2X^5+13X^4+12X^3+17X^2+14X+7)$&$7\!:\!2$\\ 
&&$\!\!\!2$&$13(13X^6+156X^5+288X^4-560X^3-936X^2+804X-148)$&\\
&&&&\\
$1185\!\!\!$&$\!\!\!A\!\!\!$&$\!\!\!1$&$(X^2-X+1)(-3X^4-3X^3+16X^2+9X-27)$&\\
&&$\!\!\!2$&$-(X-1)(X+3)(3X^4-6X^3-31X^2+34X-9)$&$2\!:\!1$\\
&&&&\\
$1197\!\!\!$&$\!\!\!A\!\!\!$&$\!\!\!1$&$(X-7)(13X+9)(X^2+X+1)(X^2-3X-3)$&\\
\end{tabular}
\newpage
\begin{tabular}{|ccclc|}\hl
&$\!\!\!C\!\!\!$&&$F(X)$&$I$\\\hl
$1325\!\!\!$&$\!\!\!A\!\!\!$&$\!\!\!1$&$5(X^6+10X^5+23X^4-10X^3+5X^2+8X+4)$&\\[6pt] 
&$\!\!\!B\!\!\!$&$\!\!\!1$&$-3(X^6+6X^5-57X^4+12X^3+171X^2+114X+13)$&\\[6pt] 
&$\!\!\!C\!\!\!$&$\!\!\!1$&$-15(X^6+6X^5-57X^4+12X^3+171X^2+114X+13)$&\\
&&&&\\
$1328\!\!\!$&$\!\!\!A^*\!\!\!$&$\!\!\!1$&$-X^6+4X^5-6X^4-2X^3-7X^2-2X-2$&\\
&&&&\\
$1360\!\!\!$&$\!\!\!A\!\!\!$&$\!\!\!1$&$-(X^2+4)(X^4+2X^3+3X^2-2X+1)$&\\
\end{tabular}
\newpage
\begin{tabular}{|ccclc|}\hl
&$\!\!\!C\!\!\!$&&$F(X)$&$I$\\\hl
$1521\!\!\!$&$\!\!\!A\!\!\!$&$\!\!\!1$&$-39(X-1)(X+3)(X^2-X+1)(X^2-X-3)$&$7\!:\!5$\\ 
&&$\!\!\!2$&$-39(X^2-X+1)(X^4+X^3+4X^2-3X+9)$&$2\!:\!1$\\ 

&&$\!\!\!3$&$-39(X-1)(X+3)(X^2+3X-1)(X^2-5X+3)$&$2\!:\!1$\\ 
&&$\!\!\!4$&$195(X^2-4X-8)(X^2+X+7)(3X^2+3X+1)$&$2\!:\!1$\\ 
&&$\!\!\!5$&$13(X-2)(5X+2)(X^2+X+7)(X^2-3X-1)$&\\
&&$\!\!\!6$&$91(3X^2+11X+9)(9X^4-39X^3+7X^2+156X+144)$&$2\!:\!5$\\
&&$\!\!\!7$&$-65(X^2-14X-3)(X^2+X-3)(23X^2-37X-9)$&$2\!:\!5$\\
\end{tabular}
\newpage
\begin{tabular}{|ccclc|}\hl
&$\!\!\!C\!\!\!$&&$F(X)$&$I$\\\hl
$1773\!\!\!$&$\!\!\!A\!\!\!$&$\!\!\!1$&$-3(81X^6-144X^5-100X^4+200X^3+80X^2-84X-36)$&\\ 
&&&&\\							 

%
$1813\!\!\!$&$\!\!\!A^*\!\!\!$&$\!\!\!1$&$7(43X^6-150X^5+261X^4-512X^3+441X^2+102X+263)$&\\
&&&&\\
$1815\!\!\!$&$\!\!\!A\!\!\!$&$\!\!\!1$&$-11X(X+4)(X^2+X-1)(X^2+X+3)$&\\ 
&&$\!\!\!2$&$-11(X^2-X+1)(X^4+X^3+8X^2-3X+9)$&$2\!:\!1$\\ 
&&$\!\!\!3$&$33(X-2)(X+2)(X^2-3X+5)(5X^2-3X+1)$&$2\!:\!1$\\
\end{tabular}
\newpage
\begin{tabular}{|ccclc|}\hl
&$\!\!\!C\!\!\!$&&$F(X)$&$I$\\\hl
$1911\!\!\!$&$\!\!\!A\!\!\!$&$\!\!\!1$&$-7(X-1)(X+3)(X^2-X+1)(X^2-X-3)$&$7\!:\!5$\\ 
&&$\!\!\!2$&$-7(X^2-X+1)(X^4+X^3+4X^2-3X+9)$&$2\!:\!1$\\ 

&&$\!\!\!3$&$-7(X-1)(X+3)(X^2+3X-1)(X^2-5X+3)$&$2\!:\!1$\\ 
&&$\!\!\!4$&$35(X^2-4X-8)(X^2+X+7)(3X^2+3X+1)$&$2\!:\!1$\\ 
&&$\!\!\!5$&$21(X-2)(5X+2)(X^2+X+7)(X^2-3X-1)$&\\
&&$\!\!\!6$&$3(3X^2+11X+9)(9X^4-39X^3+7X^2+156X+144)$&$2\!:\!5$\\
&&$\!\!\!7$&$-105(X^2-14X-3)(X^2+X-3)(23X^2-37X-9)$&$2\!:\!5$\\
\end{tabular}
\end{center}
\newpage

Peter R. Bending,\\  
Institute of Mathematics and Statistics.\\
University of Kent at Canterbury,\\
Canterbury,\\
Kent,\\
England.\\\\
Email: P.R.Bending@ukc.ac.uk 
\end{document}